\def\yes{\if00}
\def\no{\if01}
\def\iftwelvept{\yes}
\def\ifusepdf{\no}
\def\ifpsfont{\no}

\iftwelvept
\documentclass[leqno,12pt]{amsart}
\else
\documentclass[leqno]{amsart}
\fi
\usepackage{amssymb}
\usepackage{amscd}
\usepackage{latexsym}
\usepackage{verbatim}
\usepackage[all]{xy}

\ifusepdf
\usepackage{hyperref}
\else\fi
\ifpsfont
\usepackage[T1]{fontenc}
\usepackage{times}
\else\fi

\iftwelvept
\else\fi

\theoremstyle{plain}
\newtheorem{Theorem}{Theorem}[section]
\newtheorem{Proposition}[Theorem]{Proposition}
\newtheorem{Lemma}[Theorem]{Lemma}
\newtheorem{Corollary}[Theorem]{Corollary}

\newtheorem{Claim}{Claim}[Theorem]

\theoremstyle{definition}

\newtheorem{Remark}[Theorem]{Remark}

\def\rom{\textup}
\newcommand{\ZZ}{{\mathbb{Z}}}
\newcommand{\QQ}{{\mathbb{Q}}}
\newcommand{\RR}{{\mathbb{R}}}
\newcommand{\CC}{{\mathbb{C}}}
\newcommand{\PP}{{\mathbb{P}}}

\newcommand{\OO}{{\mathcal{O}}}
\newcommand{\abmod}{{\mathbb{A}}}

\newcommand{\hodge}{\lambda}
\newcommand{\aCH}{\widehat{\operatorname{CH}}}

\newcommand{\Spec}{\operatorname{Spec}}
\newcommand{\Div}{\operatorname{Div}}
\newcommand{\Supp}{\operatorname{Supp}}
\newcommand{\Hom}{\operatorname{Hom}}

\newcommand{\Ker}{\operatorname{Ker}}
\newcommand{\Coker}{\operatorname{Coker}}
\newcommand{\acherncl}{\widehat{{c}}}
\newcommand{\codim}{\operatorname{codim}}
\newcommand{\adeg}{\widehat{\operatorname{deg}}}
\newcommand{\zeros}{\operatorname{div}}
\newcommand{\Gal}{\operatorname{Gal}}
\newcommand{\trdeg}{\operatorname{tr.deg}}
\newcommand{\rank}{\operatorname{rk}}
\newcommand{\length}{\operatorname{length}}
\newcommand{\crest}{{\,\vert\,}}
\newcommand{\Aut}{\operatorname{Aut}}
\newcommand{\Sing}{\operatorname{Sing}}
\newcommand{\Proof}{{\sl Proof.}\quad}
\newcommand{\QED}{{\unskip\nobreak\hfil\penalty50\quad\null\nobreak\hfil
{$\Box$}\parfillskip0pt\finalhyphendemerits0\par\medskip}}
\newcommand{\rest}[2]{\left.{#1}\right\vert_{{#2}}}

\begin{document}

\title[The modular height of an abelian variety]%
{The modular height of an abelian variety \\
and its finiteness property}
\author{Atsushi Moriwaki}
\address{Department of Mathematics, Faculty of Science,
Kyoto University, Kyoto, 606-8502, Japan}
\email{moriwaki@math.kyoto-u.ac.jp}
\date{25/September/2003, 14:00(JP), (Version 1.0)}
\begin{abstract}
In this note, we propose the modular height
of an abelian variety defined over a field of finite
type over $\QQ$. Moreover, we prove its finiteness property.
\end{abstract}


\maketitle


\renewcommand{\theTheorem}{\arabic{Theorem}}
\renewcommand{\theClaim}{\arabic{Theorem}.\arabic{Claim}}
\renewcommand{\theequation}{\arabic{Theorem}.\arabic{Claim}}

\section*{Introduction}
In the first proof of Mordell conjecture due to Faltings,
the modular height of an abelian variety plays a crucial role.
Especially, the finiteness property of the modular height
is one of core parts for its proof.
Almost every results over a number field (i.e.,
Tate conjecture, Shafarevich conjecture, Mordell conjecture and etc)
has been generalized to a field of finite type over $\QQ$.
In this note, we propose the modular height
of an abelian variety in general and prove its finiteness property.

Let $K$ be a field of finite type over $\QQ$.
According to the paper \cite{MoArht},
we need to fix a polarization of $K$ in order to proceed with a theory of
height functions over $K$, where
a polarization of $K$ is a pair
$(B; \overline{H}_1, \ldots, \overline{H}_d)$ of
a normal and integral projective scheme $B$ over $\ZZ$ and
a sequence of nef $C^{\infty}$-hermitian line bundles 
$\overline{H}_1, \ldots, \overline{H}_d$ on $B$
with the local ring of $B$ at the generic point
isomorphic to $K$. Here we assume that $B$ is generically smooth, i.e.,
$B \times_{\ZZ} \Spec(\QQ) \to \Spec(\QQ)$ is smooth.

Let $A$ be an abelian variety over $K$. Then, using
a N\'eron model of $A$ over $B$ in codimension one,
we can introduce the Hodge sheaf $\hodge(A/K;B)$ of $A$ which is
a reflexive sheaf of rank one on $B$.
Moreover, we can give a locally integrable 
hermitian metric $\Vert\cdot\Vert_{\rm Fal}$
of $\hodge(A/K;B)$ arising from the Faltings' metric of
the good reduction part of the N\'eron model of $A$.
Then, 
\[
 \acherncl_1(\hodge(A/K;B), \Vert\cdot\Vert_{\rm Fal})
\]
can be represented by
a pair of a Weil divisor and a locally integrable function.
Thus, we can define the modular height of $A$ by the following formula:
\[
 h(A) = \adeg\left( \acherncl_1(\overline{H}_1) \cdots
 \acherncl_1(\overline{H}_d) \cdot \acherncl_1(\hodge(A/K;B), \Vert\cdot\Vert_{\rm Fal}) \right).
\]
The main purpose of this note is to prove the following result
(cf. Theorem~\ref{thm:finite:prin:abel:bound:height}):

\medskip
\begin{quotation}
If $\overline{H}_1, \ldots, \overline{H}_d$ are big, then,
for a fixed real number $c$, the set of isomorphism classes of
abelian varieties $A$ over $K$ with $h(A) \leq c$ is finite.
\end{quotation}

\medskip\noindent
Moret-Bailly \cite{Moret} proved the geometric version of the above result
using geometric intersection theory instead of Arakelov geometry.
In this sense, the above is an arithmetic generalization of
his result.

\renewcommand{\theTheorem}{\arabic{section}.\arabic{subsection}.\arabic{Theorem}}
\renewcommand{\theClaim}{\arabic{section}.\arabic{subsection}.\arabic{Theorem}.\arabic{Claim}}
\renewcommand{\theequation}{\arabic{section}.\arabic{subsection}.\arabic{Theorem}.\arabic{Claim}}

\section{Preliminaries}
\subsection{Locally integrable hermitian metric}
\setcounter{Theorem}{0}
Let $M$ be a complex manifold and $L$ a line bundle on $M$.
Let $\Vert\cdot\Vert$ be a hermitian metric of $L$, that is,
a collection of hermitian metrics of the stalks $L_x$ at all $x \in X$.
We say $\Vert\cdot\Vert$ is a {\em locally
integrable hermitian metric} (or {\em $L^1_{\rm loc}$-hermitian metric}) 
if, for each $x \in M$ and a local basis $\omega_x$ around $x$,
$\log \Vert\omega_x\Vert$ is locally integrable around $x$.
In other words, if $\Vert\cdot\Vert_0$ is a $C^{\infty}$-hermitian metric
of $L$, then $\log(\Vert\cdot\Vert/\Vert\cdot\Vert_0)$ is a locally
integrable function on $M$.

\begin{Lemma}
\label{lem:loc:int:equiv:mero}
Let $M$ be a complex manifold and $(L, \Vert\cdot\Vert)$ a hermitian line
bundle on $M$.
Let $s$ be a non-zero meromorphic section of $L$ over $M$.
Then, the hermitian metric $\Vert\cdot\Vert$ is locally integrable
if and only if so is $\log \Vert s \Vert$.
\end{Lemma}

\Proof
Let $\Vert\cdot\Vert_0$ be a $C^{\infty}$-hermitian metric of $L$.
Then,
\[
 \log \Vert s \Vert = \log(\Vert\cdot\Vert/\Vert\cdot\Vert_0) +
 \log \Vert s \Vert_0.
\]
Note that $\log \Vert s \Vert_0$ is locally integrable.
Thus, $\log \Vert s \Vert$ is locally integrable if and only if
so is $\log (\Vert\cdot\Vert/\Vert\cdot\Vert_0)$.
\QED

\begin{Lemma}
\label{lem:gen:finite:loc:int}
Let $f : Y \to X$ be a surjective, proper and generically finite
morphism of non-singular varieties over $\CC$.
Let $(L, \Vert\cdot\Vert)$ be a hermitian line bundle on $X$.
Assume that there are a non-empty 
Zariski open set $U$ of $X$ and a hermitian line
bundle $(L', \Vert\cdot\Vert')$ on $Y$ such that
$(L', \Vert\cdot\Vert')$ is isometric to
$f^*(L, \Vert\cdot\Vert)$ over $f^{-1}(U)$.
If $\Vert\cdot\Vert'$ is locally integrable, then so is
$\Vert\cdot\Vert$.
\end{Lemma}

\Proof
Shrinking $U$ if necessarily, we may assume that $f$ is \'etale over $U$.
We set $V = f^{-1}(U)$.
Let $s$ be a non-zero rational section of $L$.
Note that there is a divisor $D$ on $Y$ such that
$L' = f^*(L) \otimes \OO_Y(D)$ and $\Supp(D) \subseteq Y \setminus V$.
Thus, $f^*(s)$ gives rise to a rational section $s'$ of $L'$.
Then, $\log \Vert s' \Vert'$ is locally integrable 
by Lemma~\ref{lem:loc:int:equiv:mero}.
Since $\rest{f^*(\log \Vert s \Vert)}{V} = 
\rest{\log \Vert s' \Vert'}{V}$, we can see that
$f^*(\log \Vert s \Vert)$ is locally integrable.
Let $[f^*(\log \Vert s \Vert)]$ be a current associated to
the locally integrable function $f^*(\log \Vert s \Vert)$.
Then,
\cite[Proposition~1.2.5]{KMSemi},
there is a locally integrable function $g$ on $X$ with
$f_*[f^*(\log \Vert s \Vert)] = [g]$.
Since $f$ is \'etale over $U$,
we can easily see that
\[
 (\rest{f}{V})_*[(\rest{f}{V})^*(\rest{\log \Vert s \Vert}{U})] =
 \deg(f)[\rest{\log \Vert s \Vert}{U}].
\]
Thus, $g = \deg(f) \log \Vert s \Vert$ almost everywhere over $U$.
Therefore, so is over $X$ because $U$ is a non-empty Zariski open set of $X$.
Hence, $\log \Vert s \Vert$ is locally integrable on $X$.
\QED

\subsection{Hermitian metric with logarithmic singularities}
\setcounter{Theorem}{0}
Let $X$ be a normal variety over $\CC$ and $Y$ a proper closed subscheme of $X$.
Let $(L, \Vert\cdot\Vert)$ be a hermitian line bundle on $X$. We say 
$(L, \Vert\cdot\Vert)$ is a {\em $C^{\infty}$-hermitian line bundle 
with logarithmic singularities along $Y$} 
if the following conditions are satisfied:
\begin{enumerate}
\renewcommand{\labelenumi}{(\arabic{enumi})}
\item 
$\Vert\cdot\Vert$ is $C^{\infty}$ over $X \setminus Y$.

\item
Let $\Vert\cdot\Vert_0$ be a $C^{\infty}$-hermitian metric of $L$.
For each $x \in Y$, let $f_1, \ldots, f_m$ be a system of local equations of $Y$ 
around $x$, i.e.,
$Y$ is given by $\{ z \in X \mid f_1(z) = \cdots = f_m(z) = 0\}$ around $x$. 
Then, there are positive constants $C$ and $r$ such that
\[
 \max \left\{ \frac{\Vert\cdot\Vert}{\Vert\cdot\Vert_0}, 
         \frac{\Vert\cdot\Vert_0}{\Vert\cdot\Vert} \right\}
 \leq C \left(-\sum_{i=1}^m \log \vert f_i \vert\right)^r
\]
around $x$.
\end{enumerate}
Note that the above definition does not depend on the choice
of the system of local equations
$f_1, \ldots, f_m$. Moreover,
it is easy to see that if $(L, \Vert\cdot\Vert)$ is a 
$C^{\infty}$-hermitian line bundle with logarithmic singularities along $Y$,
then $\Vert\cdot\Vert$ is locally integrable.

\begin{Lemma}
\label{lem:log:sing:herm:pullback}
Let $\pi : X' \to X$ be a proper morphism
of normal varieties over $\CC$ and $Y$ a proper closed subscheme of $X$.
Let $(L, \Vert\cdot\Vert)$ be a hermitian line bundle on $X$ such that
$\Vert\cdot\Vert$ is $C^{\infty}$ over $X \setminus Y$.
If $\pi(X') \not\subseteq Y$ and $(L, \Vert\cdot\Vert)$ 
has logarithmic singularities along $Y$, 
then so does $\pi^*(L, \Vert\cdot\Vert)$
along $\pi^{-1}(Y)$.
Moreover, if $\pi$ is surjective and $\pi^*(L, \Vert\cdot\Vert)$ has
logarithmic singularities along $\pi^{-1}(Y)$, then
so does $(L, \Vert\cdot\Vert)$ along $Y$.
\end{Lemma}

\Proof
Let $\{ f_1, \ldots, f_m\}$ be a system of local equations of $Y$.
Then, $\{ \pi^*(f_1), \ldots, \pi^*(f_m)\}$ is a system of local equation
of $\pi^{-1}(Y)$. Thus, our assertion is obvious.
\QED

\subsection{Faltings' metric}
\setcounter{Theorem}{0}
Let $X$ be a normal variety over $\CC$.
Let $f : A \to X$ be a $g$-dimensional semi-abelian scheme over $X$.
We assume that there is a non-empty Zariski open set $U$ of $X$ such that
$f$ is an abelian scheme over $U$.
Let $\hodge_{A/X}$ be the Hodge line bundle of $A \to X$, i.e.,
\[
 \lambda_{A/X} = \det \left( \epsilon^* \left( \Omega_{A/X} \right) \right),
\]
where $\epsilon : X \to A$ is the identity of the semi-abelian scheme 
$A \to X$.
At each $x \in U$, we can give a hermitian metric of
$(\hodge_{A/X})_x$ in the following way:
For $\alpha \in \bigwedge^g H^0(\Omega_{A_x})$,
\[
 (\Vert \alpha \Vert_x)^2 = \left( \frac{\sqrt{-1}}{2} \right)^g
 \int_{A_x} \alpha \wedge \bar{\alpha}.
\]
Then, a collection of metrics $\{ \Vert\cdot\Vert_x \}_{x \in U}$
gives rise to a $C^{\infty}$-hermitian metrics $\Vert\cdot\Vert_{\rm Fal}$
of $\rest{\hodge_{A/X}}{U}$.
Moreover, it is well-known that $\Vert\cdot\Vert_{\rm Fal}$ extends to
a $C^{\infty}$-hermitian metric of $\hodge_{A/X}$ with logarithmic singularities
along $X \setminus U$ (cf. \cite[Th\'eor\`em~3.2 in Expos\'e~I]{Mordell}).
By abuse of notation, this extended metric is also denoted by 
$\Vert\cdot\Vert_{\rm Fal}$
and is called Faltings' metric of $\hodge_{A/X}$.

\begin{Lemma}
\label{lem:loc:L1:Faltings:metric}
Let $X$ be a smooth variety over $\CC$ and $X_0$
a non-empty Zariski open set of $X$.
Let $A_0 \to X_0$ be an abelian scheme over $X_0$.
Let $\lambda$ be a line bundle on $X$ such that
$\rest{\lambda}{X_0}$ gives rise to the Hodge line bundle $\hodge_{A_0/X_0}$
of $A_0 \to X_0$.
Then, Faltings' metric $\Vert\cdot\Vert_{\rm Fal}$ of $\hodge_{A_0/X_0}$
over $X_0$
extends to a locally integrable metric of $\hodge$ over $X$.
\end{Lemma}

\Proof
By virtue of Lemma~\ref{lem:Gabber} (Gabber's lemma),
there is a proper, surjective and generically finite morphism
$\pi : X' \to X$ of smooth varieties over $\CC$ such that
the abelian scheme $A_0 \times_{X_0} \pi^{-1}(X_0)$ over $\pi^{-1}(X_0)$
extends to a semi-abelian scheme $f' : A' \to X'$.
Let $\hodge_{A'/X'}$ be the Hodge line bundle of $A' \to X'$ and
$\Vert\cdot\Vert'_{\rm Fal}$ Faltings' metric of $\hodge_{A'/X'}$.
Then, $\rest{(\hodge_{A'/X'}, \Vert\cdot\Vert'_{\rm Fal})}{X'_0}$ is
isometric to $\pi_0^*
(\hodge_{A_0/X_0}, \Vert\cdot\Vert_{\rm Fal})$, 
where $X'_0 = \pi^{-1}(X_0)$ and $\pi_0 = \rest{\pi}{X'_0}$.
Therefore, by Lemma~\ref{lem:gen:finite:loc:int},
$\Vert\cdot\Vert_{\rm Fal}$ extends to a locally integrable metric over $X$.
\QED

\subsection{N\'eron model}
\setcounter{Theorem}{0}
Let $R$ be a discrete valuation ring and $K$ the quotient field of $R$.
Let $A$ be an abelian variety over $K$. Then, there is
a smooth group scheme $\mathcal{A} \to \Spec(R)$ 
of finite type over $R$
with the following properties (cf. \cite{Neron}):
\begin{enumerate}
\renewcommand{\labelenumi}{(\arabic{enumi})}
\item
The generic fiber of $\mathcal{A} \to \Spec(R)$ is $A$.

\item
Let $\mathcal{X} \to \Spec(R)$ be a smooth scheme over $R$ and
$X$ the generic fiber of $\mathcal{X} \to \Spec(R)$.
Then, any morphism $X \to A$ over $K$ extends uniquely to
a morphism $\mathcal{X} \to \mathcal{A}$ over $R$.
\end{enumerate}
The smooth group scheme $\mathcal{A} \to \Spec(R)$ is called
the N\'eron model of $A$ over $R$.
We would like to generalize it to a higher dimensional base scheme.

\medskip
Let $B$ be an irreducible noetherian normal scheme and
$K$ the function field of $B$, i.e.,
the local ring at the generic point of $B$.
Let $A$ be an abelian variety over $K$.
A smooth group scheme $f : \mathcal{A} \to B$
is called the N\'eron model of $A$ over $B$ if
(1) $f : \mathcal{A} \to B$ is of finite type over $B$
and (2) for every point $x \in B$ of codimension one,
$\rest{\mathcal{A}}{\Spec(\OO_x)} \to \Spec(\OO_x)$
is the N\'eron model of $A$ over $\Spec(\OO_x)$.
Let $\mathcal{X} \to B$ a smooth scheme over $B$ and
$X$ the generic fiber of $\mathcal{X} \to B$.
Let $\phi_K : X \to A$ be a morphism over $K$.
If $f : \mathcal{A} \to B$ is the N\'eron model of $A$,
then, by the property (2) and Weil's extension theorem
(cf. \cite[Theorem~1 in 4.4]{Neron}),
there is the unique extension $\phi : \mathcal{X} \to \mathcal{A}$
of $\phi_K$ over $B$.

\begin{Proposition}
\label{prop:existence:Neron:model}
Let $B$, $K$ and $A$ be same as above.
Then there is a non-empty big open set $B'$ of $B$
\rom{(}i.e., $\codim(B \setminus B') \geq 2$\rom{)} such that
a N\'eron model of $A$ over $B'$ exists.
This N\'eron model is called a N\'eron model of $A$
over $B$  in codimension one.
\end{Proposition}

\Proof
First of all, we can take a non-empty Zariski open set
$B_0$ of $B$ and an abelian scheme $\mathcal{A}_0 \to B_0$
whose generic fiber is $A$.
Let
\[
 B \setminus B_0 = D_1 \cup D_2 \cup \cdots \cup D_n
\]
be the irreducible decomposition of $B \setminus B_0$.
We assume that $\codim(D_i) = 1$ for $1 \leq i \leq r$ and
$\codim(D_j) \geq 2$ for $r < j \leq n$.
Let $x_i$ be the generic point of $D_i$.
For $1 \leq i \leq r$,
let $A_i \to \Spec(\OO_{x_i})$ be the N\'eron model of $A$
over $\Spec(\OO_{x_i})$.
Then, there are an open set $B_i$ containing $x_i$ and
a smooth group scheme $\mathcal{A}_i \to B_i$
as the extension of $A_i \to \Spec(\OO_{x_i})$
such that $\mathcal{A}_i$ is of finite type over $B_i$.
Replacing $B_i$ by $B_i \setminus (D_1 \cup \cdots
\cup D_{i-1} \cup D_{i+1} \cup \cdots \cup D_n)$,
we may assume that
\[
 B_i \cap (D_1 \cup \cdots
\cup D_{i-1} \cup D_{i+1} \cup \cdots \cup D_n) = \emptyset.
\]
By \cite[Lemma~3.3 in Chapter~I]{FalRat} or 
\cite[Proposition~2.7 in Chapter~1]{FalChai},
$\mathcal{A}_0 \to B_0$ coincides with $\mathcal{A}_i \to B_i$
over $B_0 \cap B_i$. Moreover,
$B_i \cap B_j \subseteq B_0$
for $1 \leq i \not= j \leq r$.
Therefore, if we set $B' = B_0 \cup B_1 \cup \cdots \cup B_r$,
we have our desired smooth group scheme $\mathcal{A} \to B'$.
\QED

\subsection{Semi-abelian reduction}
\setcounter{Theorem}{0}
Let $B$ be an irreducible normal noetherian scheme and
$K$ the local ring at the generic point of $B$.
Let $A$ be an abelian variety over $K$.
We say $A$ has semi-abelian reduction over $B$ in codimension one
if there are a big open set $B_1$ of $B$
(i.e., $\codim(B \setminus B') \geq 2$) and
a semi-abelian scheme $\mathcal{A} \to B_1$ such that
the generic fiber of $\mathcal{A} \to B_1$ is $A$.

\begin{Proposition}
\label{prop:semi:abelian:reduction}
Let $B$, $K$ and $A$ be same as above.
Let $m$ be a positive integer which has a factorization
$m = m_1m_2$ with $m_1, m_2 \geq 3$ and $m_1$ and $m_2$ relatively prime
\rom{(}for example $m=12 = 3 \cdot 4$\rom{)}.
If $A[m](\overline{K}) \subseteq A(K)$,  then
$A$ has semi-abelian reduction in codimension one over $B$.
\end{Proposition}

\Proof
Let $x$ be a codimension one point of $B$.
Then, there is $m_i$ which is not divisible by the characteristic of 
the residue field of $\OO_{B,x}$.
Moreover, $A[m_i](\overline{K}) \subseteq A(K)$.
Thus, by \cite[expos\'e 1, Corollaire~5.18]{Spgenus},
$A$ has semi-abelian reduction
at $x$.

Let $B_0$ be a non-empty Zariski open set
of $B$ such that we can take an abelian scheme $\mathcal{A}_0 \to B_0$
whose generic fiber is $A$.
Let
\[
 B \setminus B_0 = D_1 \cup D_2 \cup \cdots \cup D_n
\]
be the irreducible decomposition of $B \setminus B_0$.
We assume that $\codim(D_i) = 1$ for $1 \leq i \leq r$ and
$\codim(D_j) \geq 2$ for $r < j \leq n$.
Let $x_i$ be the generic point of $D_i$.
Then, for each $i=1, \ldots, r$,
there are an open set $B_i$ of $B$ and
a semi-abelian scheme $\mathcal{A}_i \to B_i$ with
$x_i \in B_i$.
Shrinking $B_i$ if necessarily,
we may assume that
$B_i \cap B_j \subseteq B_0$ for all $1 \leq i \not= j \leq r$.
Thus, as in Proposition~\ref{prop:existence:Neron:model},
if we set $B' = B_0 \cup B_1 \cup \cdots \cup B_r$,
then we have our desired semi-abelian scheme
$\mathcal{A} \to B'$.
\QED

\begin{Lemma}[Gabber's lemma]
\label{lem:Gabber}
Let $U$ be a dense Zariski open set of an integral, normal and excellent scheme 
$S$ and $A$ an abelian scheme over $U$. Then, there is
a proper, surjective and generically finite morphism
$\pi : S' \to S$ of integral, normal and excellent schemes 
such that the abelian scheme $A \times_U f^{-1}(U)$ over $f^{-1}(U)$
extends to an semi-abelian scheme over $S'$
\end{Lemma}

\Proof
In \cite[Th\'eor\`em and Proposition~4.10 in Expos\'e~V]{Mordell}, 
the existence of $\pi : S' \to S$ and the extension of the abelian scheme
is proved under the assumption
$\pi : S' \to S$ is proper and surjective.
Let $S'_{\eta}$ be the generic fiber of $\pi$.
Let $z$ be the closed point of $S'_{\eta}$ and
$Z$ the closure of $z$ in $S'$.
Moreover, let $S_1$ be the normalization of $Z$.
Then, $\pi_1 : S_1 \to Z \to S$ is our desired morphism.
\QED

\subsection{The Hodge sheaf of an abelian variety}
\label{subsec:mod:sheaf:ab:var}
Let $G \to S$ be a smooth group scheme over $S$.
Then, {\em the Hodge line bundle $\hodge_{G/S}$ of $G \to S$} is 
given by
\[
 \hodge_{G/S} = \det \left( \epsilon^*\left( \Omega_{G/S} \right) \right),
\]
where $\epsilon$ is the identity of the group scheme
$G \to S$.

Let $B$ be an irreducible and normal noetherian scheme.
Let $K$ be the function field of $B$ 
(i.e., the local ring at the generic point).
Let $A$ be an abelian variety over $K$.
By Proposition~\ref{prop:existence:Neron:model},
there is a big open set $B'$ of $B$ such that
the N\'eron model $\mathcal{A}' \to B'$ of $A$ over $B'$ exists.
Let $\iota : B' \to B$ be the natural inclusion map.
{\em The Hodge sheaf $\hodge(A/K; B)$ of $A$ with respect to $B$}
is defined by
\[
\hodge(A/K; B) = 
\iota_*\left( \hodge_{\mathcal{A}'/B'} \right).
\]
Note that $\hodge(A/K; B)$ is a reflexive sheaf of rank one on $B$.

From now on, we assume that
the characteristic of $K$ is zero.
Let $\phi : A \to A'$ be an isogeny of abelian varieties
over $K$. Let $\mathcal{A}$ and $\mathcal{A}'$ be the
N\'eron models in codimension one over $B$ of $A$ and
$A'$ respectively.
Since there is an injective homomorphism
\[
\phi^* : \hodge(A'/K;B) \to \hodge(A/K;B),
\]
we can find an effective Weil divisor $D_{\phi}$ such that
\[
c_1(\hodge(A'/K;B)) + D_{\phi} = c_1(\hodge(A/K;B)).
\]
The ideal sheaf $\OO_B(-D_{\phi})$ is denoted by
$\mathcal{I}_{\phi}$.

\begin{Lemma}
\label{lemma:mult:iso:isodual:deg}
Let $\phi^{\vee} : {A'}^{\vee} \to A^{\vee}$ be the dual of
$\phi : A \to A'$.
We assume that $B = \Spec(R)$ for some discrete valuation ring $R$
and that
$A$ and $A'$ have semi-abelian reduction over $R$.
Then, $\mathcal{I}_{\phi} \cdot \mathcal{I}_{\phi^{\vee}} =
\deg(\phi) R$.
\end{Lemma}

\Proof
Let $R'$ be an extension of $R$ such that
$R'$ is a complete discrete valuation ring and
the residue field of $R'$ is algebraically closed.
Then, by \cite[Expos\'e~VII, Th\'ero\`em~2.1.1]{Mordell},
$(\mathcal{I}_{\phi} \cdot \mathcal{I}_{\phi^{\vee}}) R' =
\deg(\phi) R'$. Here $R'$ is faithfully flat over $R$. Thus,
$\mathcal{I}_{\phi} \cdot \mathcal{I}_{\phi^{\vee}} =
\deg(\phi) R$.
\QED

\subsection{The moduli of abelian varieties}
\setcounter{Theorem}{0}
To prove the finiteness property of the modular height,
it is very important to get a good compactification of the moduli
space of abelian varieties. 
For simplicity, an abelian variety with a polarization of degree $l^2$
is called an {\em $l$-polarized abelian variety}.

\begin{Theorem}
\label{thm:universal:extension:abel:semiable:with:level}
Let $g$, $l$ and $m$ be positive integers with $m \geq 3$.
Let $\abmod_{g,l,m,\QQ}$ be the moduli space of $g$-dimensional 
and $l$-polarized abelian
varieties over $\QQ$ with an
$m$-level structure.
Then, there exists \rom{(a)} normal projective arithmetic varieties
$\abmod_{g,l,m}^*$ and $Y^*$ \rom{(}i.e.,
$\abmod_{g,l,m}^*$ and $Y^*$ are normal and integral schemes flat and
projective over $\ZZ$\rom{)}, \rom{(b)} a surjective and
generically finite morphism
$f : Y^* \to \abmod_{g,l,m}^*$, \rom{(c)} a positive integer $n$,
\rom{(d)} a line bundle $L$ on $\abmod_{g,l,m}^*$, and 
\rom{(e)} a semi-abelian scheme $G \to Y^*$
with the following properties:
\begin{enumerate}
\renewcommand{\labelenumi}{(\arabic{enumi})}
\item
$\abmod_{g,l,m,\QQ}$ is a Zariski open set of 
$\abmod_{g,l,m,\QQ}^* = \abmod_{g,l,m}^* \times_{\ZZ} \Spec(\QQ)$ and
$L$ is very ample on $\abmod_{g,l,m}^*$.

\item
Let $\lambda_{G/Y^*}$ be the Hodge line bundle of the semi-abelian scheme
$G \to Y^*$. Then,
$f^*(L) = \hodge_{G/Y^*}^{\otimes n}$ on 
$Y^*_{\QQ}= Y^* \times_{\ZZ} \Spec(\QQ)$.

\item
Let $U_{\QQ} \to \abmod_{g,l,m,\QQ}$ be the universal 
$g$-dimensional and $l$-principally polarized
abelian scheme with an $m$-level structure.
Let $Y_{\QQ}$ be the pull-back of $\abmod_{g,l,m,\QQ}$ by
$f_{\QQ} : Y^*_{\QQ} \to \abmod_{g,l,m,\QQ}^*$, i.e.,
$Y_{\QQ} = (f_{\QQ})^{-1}(\abmod_{g,l,m,\QQ})$.
Then, $G_{\QQ} \to Y^*_{\QQ}$ is an extension of
the abelian scheme $U_{\QQ} \times_{\abmod_{g,l,m,\QQ}} Y_{\QQ} \to Y_{\QQ}$.
\rom{(}Note that $\rest{G}{Y_{\QQ}} \to Y_{\QQ}$ is naturally a
$g$-dimensional and $l$-polarized
abelian scheme with
an $m$-level structure.\rom{)}

\item
$L$ has a metric $\Vert\cdot\Vert$ over $\abmod_{g,l,m,\QQ}(\CC)$
such that $f^*((L, \Vert\cdot\Vert))$ is isometric to  
$\left(\hodge_{G/Y^*}, \Vert\cdot\Vert_{\rm Fal}
\right)^{\otimes n}$
over $Y_{\QQ}(\CC)$.
Moreover, $\Vert\cdot\Vert$ has logarithmic singularities along
$\abmod_{g,l,m,\QQ}^*(\CC) \setminus \abmod_{g,l,m,\QQ}(\CC)$.
\end{enumerate}
\end{Theorem}

\Proof
Let $U_{\QQ} \to \abmod_{g,l,m,\QQ}$ be the universal $l$-polarized
abelian scheme with an $m$-level structure.
By \cite[Th\'eor\`eme~2.2 in Expos\'e~IV]{Mordell},
there are a normal projective variety $\abmod_{g,l,m,\QQ}^*$,
a positive integer $n$ and
a very ample line bundle $L_{\QQ}$ on $\abmod_{g,l,m,\QQ}^*$ with
the following properties:
\begin{enumerate}
\renewcommand{\labelenumi}{(\roman{enumi})}
\item
$\abmod_{g,l,m,\QQ}$ is an Zariski open set of
$\abmod_{g,l,m,\QQ}^*$.

\item
By Gabber's lemma (cf. Lemma~\ref{lem:Gabber}),
there is a surjective and generically finite morphism
$h_{\QQ} : S'_{\QQ} \to \abmod_{g,l,m,\QQ}^*$ of normal
projective varieties over $\QQ$ such that
the abelian scheme
$U_{\QQ} \times_{\abmod_{g,l,m,\QQ}} 
h_{\QQ}^{-1}(\abmod_{g,l,m,\QQ}) \to h_{\QQ}^{-1}(\abmod_{g,l,m,\QQ})$
extends to a semi-abelian scheme $G'_{\QQ} \to S'_{\QQ}$.
Then, $h_{\QQ}^*(L_{\QQ}) = \hodge_{G'_{\QQ}/S'_{\QQ}}^{\otimes n}$.
\end{enumerate}
Since $L_{\QQ}$ is very ample, there is an embedding
$\abmod_{g,l,m,\QQ}^* \hookrightarrow \PP^N_{\QQ}$ in terms of
$L_{\QQ}$.
Let $\abmod_{g,l,m}^*$ be the closure of the image of
\[
\abmod_{g,l,m,\QQ}^* \hookrightarrow \PP^N_{\QQ} \to \PP^N_{\ZZ}.
\]
Moreover, let $L$ be the pull-back of $\OO_{\PP^N_{\ZZ}}(1)$ by
the embedding $\abmod_{g,l,m}^* \hookrightarrow \PP^N_{\ZZ}$.
Then, $\abmod_{g,l,m,\QQ}^* = \abmod_{g,l,m}^* \times_{\ZZ} \Spec(\QQ)$ and
$L_{\QQ} = \rest{L}{\abmod_{g,l,m,\QQ}^*}$.
Let $S'$ be the normalization of $\abmod_{g,l,m}^*$ in the function field of $S'_{\QQ}$.
Then, there is an open set $S'_0$ of $S'$ such that
$G'$ is an ableian scheme over $S'_0$ and $G' \times_{S'} S'_0 \to
S'_0$ coincides with
the abelian scheme $U_{\QQ} \times_{\abmod_{g,l,m,\QQ}} 
h_{\QQ}^{-1}(\abmod_{g,l,m,\QQ}) \to h_{\QQ}^{-1}(\abmod_{g,l,m,\QQ})$
over $\QQ$.
Thus, using Gabber's lemma again,
there are a surjective and generically finite 
morphism of normal arithmetic varieties
$h_2 : Y^* \to S'$ and a semi-abelian scheme $G \to Y^*$
such that $G \to Y^*$ is an extension of $G' \times_{S'} 
h_2^{-1}(S'_0) \to h_2^{-1}(S'_0)$. We set $Y^*_{\QQ} = 
Y^* \times_{\ZZ} \Spec(\QQ)$. Then,
$G$ over $Y^*_{\QQ}$ 
is equal to
$G'_{\QQ} \times_{ S'_{\QQ} } Y^*_{\QQ} \to 
Y^*_{\QQ}$ 
by the uniqueness of
semi-abelian extension.
Thus, if we set $f = h \cdot h_1$, then  $f^*(L) = \hodge_{G/Y^*}^{\otimes n}$
over $Y^*_{\QQ}$.

Finally, since
$\rest{L_{\QQ}}{\abmod_{g,l,m,\QQ}} = 
\hodge_{U_{\QQ}/\abmod_{g,l,m,\QQ}}^{\otimes n}$,
if we give $L_{\QQ}$ a metric arising from the Faltings' metric
of $\hodge_{U_{\QQ}/\abmod_{g,l,m,\QQ}}$, then assertion of (4) follows from
Lemma~\ref{lem:log:sing:herm:pullback} and 
\cite[Th\'eor\`em~3.2 in Expos\'e~I]{Mordell}.
\QED

\subsection{Arakelov geometry}
\setcounter{Theorem}{0}
In this note, a flat and quasi-projective integral
scheme over $\ZZ$ is called an {\em arithmetic variety}.
If it is smooth over $\QQ$, then it is said to be {\em generically smooth}.

Let $X$ be a generically smooth arithmetic variety.
A pair $(Z, g)$ is called an {\em arithmetic cycle of codimension $p$}
if $Z$ is a cycle of codimension $p$ and
$g$ is a current of type $(p-1,p-1)$ on $X(\CC)$.
We denote by $\widehat{Z}^p(X)$ the set of all arithmetic cycles
on $X$. We set 
\[
\aCH^p(X) = \widehat{Z}^p(X)/\!\!\sim,
\]
where $\sim$ is the arithmetic linear equivalence.

Let $\overline{L} = (L, \Vert\cdot\Vert)$
be a $C^{\infty}$-hermitian line bundle on $X$.
Then, a homomorphism
\[
\acherncl_1(\overline{L}) \cdot :
\aCH^p(X) \to \aCH^{p+1}(X)
\]
is define by
\[
\acherncl_1(\overline{L}) \cdot (Z, g) =
\left(\text{$\zeros(s)$ on $Z$}, [-\log (\Vert s \Vert^2_Z)]
+ c_1(\overline{L}) \wedge g \right),
\]
where $s$ is a rational section of $\rest{L}{Z}$ and
$[-\log (\Vert s \Vert^2_Z)]$ is a current given by
$\phi \mapsto -\int_{Z(\CC)}\log(\Vert s \Vert^2_Z)\phi$.

Here we assume that $X$ is projective.
Then we can define the arithmetic degree map
\[
\adeg : \aCH^{\dim X}(X) \to \RR
\]
by
\[
\adeg \left( \sum_P n_P P, g \right) = \sum_{P} n_P \log(\#(\kappa(P))) +
\frac{1}{2} \int_{X(\CC)} g.
\]
Thus, if $C^{\infty}$-hermitian line bundles
$\overline{L}_1, \ldots, \overline{L}_{\dim X}$ are given, then
we can get the number
\[
\adeg \left(\acherncl_1(\overline{L}_1) \cdots 
\acherncl_1(\overline{L}_{\dim X}) \right),
\]
which is called the {\em arithmetic intersection number
of $\overline{L}_1, \ldots, \overline{L}_{\dim X}$}.

\medskip
Let $X$ be a projective arithmetic variety.
Note that $X$ is not necessarily generically smooth.
Let $\overline{L}_1, \ldots, \overline{L}_{\dim X}$ be
$C^{\infty}$-hermitian line bundles on $X$.
By \cite{Hiro},
we can find a generic resolution of singularities $\mu : Y \to X$, i.e.,
$\mu : Y \to X$ is a projective and birational morphism such that
$Y$ is a generically smooth projective
arithmetic variety.
Then, we can see that the arithmetic intersection number
\[
\adeg \left(
\acherncl_1(\mu^*(\overline{L}_1)) \cdots
\acherncl_1(\mu^*(\overline{L}_{\dim X})) \right)
\]
does not depend on the choice of the generic
resolution of singularities $\mu : Y \to X$.
Thus, we denote this number by
\[
\adeg \left(
\acherncl_1(\overline{L}_1) \cdots
\acherncl_1(\overline{L}_{\dim X}) \right).
\]

\medskip
Let $\overline{L}_1, \ldots, \overline{L}_l$ be
$C^{\infty}$-hermitian line bundles on a projective arithmetic
variety $X$.
Let $V$ be an $l$-dimensional integral closed subscheme on $X$.
Then, 
$\adeg\left(\acherncl_1(\overline{L}_1) \cdots \acherncl_1(\overline{L}_l)
\crest V\right)$ is defined by
\[
\adeg\left(\acherncl_1(\rest{\overline{L}_1}{V}) \cdots 
\acherncl_1(\rest{\overline{L}_l}{V})\right).
\]
Moreover, for an $l$-dimensional cycle
$Z = \sum_{i} n_i V_i$ on $X$,
$\adeg\left(\acherncl_1(\overline{L}_1) \cdots \acherncl_1(\overline{L}_l)
\crest Z \right)$ is defined by
\[
 \sum_i 
n_i \adeg\left(\acherncl_1(\overline{L}_1) \cdots \acherncl_1(\overline{L}_l)
\crest V_i \right).
\]

\subsection{The positivity of $C^{\infty}$-hermitian 
$\QQ$-line bundles on a projective arithmetic variety}
\setcounter{Theorem}{0}
Let $X$ be a projective arithmetic variety and
$\overline{L}$ a $C^{\infty}$-hermitian $\QQ$-line bundle on $X$.
Let us introduce several kinds of
the positivity of $C^{\infty}$-hermitian $\QQ$-line bundles.

$\bullet${\bf ample}:
We say $\overline{L}$ is {\em ample} if
$L$ is ample on $X$, $c_1(\overline{L})$ is positive form on $X(\CC)$, and
there is a positive number $n$ such that
$L^{\otimes n}$ is generated by
the set $\{ s \in H^0(X, L^{\otimes n}) \mid \Vert s \Vert_{\sup} < 1 \}$.

$\bullet${\bf nef}:
We say  $\overline{L}$ is {\em nef} if
$c_1(\overline{L})$ is a semipositive form on $X(\CC)$ and,
for all one-dimensional integral closed subschemes $\Gamma$ of $X$,
$\adeg \left( \acherncl_1(\overline{L}) \crest \Gamma \right) \geq 0$.

$\bullet${\bf big}:
$\overline{L}$ is said to be {\em big} if 
$\rank_{\ZZ} H^0(X, L^{\otimes m}) = O(m^{\dim X_{\QQ}})$ and there is a non-zero
section $s$ of $H^0(X, L^{\otimes n})$ with $\Vert s \Vert_{\sup} < 1$ for
some positive integer $n$.

$\bullet${\bf $\pmb{\QQ}$-effective}:
$\overline{L}$
is said to be {\em $\QQ$-effective} if there is a positive integer $n$ and
a non-zero $s \in H^0(X, L^{\otimes n})$ with $\Vert s \Vert_{\sup} \leq 1$.

$\bullet${\bf pseudo-effective}:
$\overline{L}$
is said to be {\em pseudo-effective} if there are
(1) a sequence $\{ \overline{L}_n\}_{n=1}^{\infty}$ of $\QQ$-effective 
$C^{\infty}$-hermitian $\QQ$-line bundles,
(2) $C^{\infty}$-hermitian $\QQ$-line bundles 
$\overline{E}_1, \ldots, \overline{E}_r$ and
(3) sequences $\{ a_{1, n} \}_{n=1}^{\infty}, \ldots, \{a_{r, n}\}_{n=1}^{\infty}$
of rational numbers such that
\[
\acherncl_1(\overline{L}) = \acherncl_1(\overline{L}_n) + 
\sum_{i=1}^r a_{i,n}\acherncl_1(\overline{E}_i)
\]
in $\aCH^1(X) \otimes \QQ$ and
$\lim_{n\to\infty} a_{i, n} = 0$ for all $i$.
If $\overline{L}_1 \otimes \overline{L}_2^{\otimes -1}$ is pseudo-effective
for $C^{\infty}$-hermitian $\QQ$-line bundles $\overline{L}_1, \overline{L}_2$
on $X$, then we denote this by $\overline{L}_1 \succsim \overline{L}_2$.

\subsection{Polarization of a finitely generated field over $\QQ$}
\label{subsec:polarization}
\setcounter{Theorem}{0}
Let $K$ be a field of finite type over the rational number field $\QQ$
with $d = \trdeg_{\QQ}(K)$.
A pair $\overline{B} = (B; \overline{H}_1, \ldots, \overline{H}_d)$
of a normal projective arithmetic variety $B$ and
a sequence $\overline{H}_1, \ldots, \overline{H}_d$ of
$C^{\infty}$-hermitian line bundles on $B$ is called a polarization if
the function field of $B$ (i.e., the local ring at the generic point)
is $K$ and
$\overline{H}_1, \ldots, \overline{H}_d$ are nef.
Here $\deg(\overline{B})$ is given by
\[
 \int_{B(\CC)} c_1(\overline{H}_1) \wedge \cdots \wedge
 c_1(\overline{H}_d).
\]
Namely,
\[
 \deg(\overline{B}) = \begin{cases}
 [K : \QQ] & \text{if $d = 0$}, \\
 \text{$\deg((H_1)_{\QQ} \cdots (H_d)_{\QQ})$ on 
 $B \times_{\ZZ} \Spec(\QQ)$} & \text{if $d > 0$}.
 \end{cases}
\]
If $B$ is generically smooth, then the polarization $\overline{B}$ is 
said to be {\em generically smooth}. Moreover, we say the polarization 
$\overline{B} = (B; \overline{H}_1, \ldots, \overline{H}_d)$ is 
{\em fine} (resp. {\em strictly fine})
if there is a generically finite morphism $\pi : B' \to B$
of normal projective arithmetic varieties, 
a generically finite morphism $\mu : B' \to (\PP^1_{\ZZ})^d$ and
ample $C^{\infty}$-hermitian $\QQ$-line bundles
$\overline{L}_1, \ldots, \overline{L}_d$ on $\PP^1_{\ZZ}$
such that 
$\pi^*(\overline{H}_i) \otimes \mu^*(p_i^*(\overline{L}_i))^{\otimes -1}$
is pseudo-effective (resp. $\QQ$-effective) 
for every $i$, where $p_i : (\PP^1_{\ZZ})^d \to \PP^1_{\ZZ}$
is the projection to the $i$-th factor.
Note that if $\overline{H}_1, \ldots, \overline{H}_d$
are big, then the polarization
$(B; \overline{H}_1, \ldots, \overline{H}_d)$ is
strictly fine. Moreover, if $\overline{B}$ is fine,
then $\deg(\overline{B}) > 0$.

Let us see the following proposition.

\begin{Proposition}
\label{prop:finite:divisor:fine:polarization}
Let $\overline{B} = (B; \overline{H}_1, \ldots, \overline{H}_d)$
be a strictly fine polarization of $K$.
Then,
for all $h$, the number of prime divisors $\Gamma$ on $B$ with
\[
\adeg(\acherncl_1(\overline{H}_1) \cdots \acherncl_1(\overline{H}_1) \crest \Gamma) \leq h
\]
is finite.
\end{Proposition}

\Proof
Let us begin with the following lemma.

\begin{Lemma}
\label{lem:equiv:div:finite}
Let $\pi : X' \to X$ be a generically finite morphism of normal projective arithmetic varieties.
Let $\overline{H}_1, \ldots, \overline{H}_d$ be nef
$C^{\infty}$-hermitian line bundles on $X$, where $d = \dim X_{\QQ}$.
Then, the following are equivalent:
\begin{enumerate}
\renewcommand{\labelenumi}{(\arabic{enumi})}
\item
For all $h$, the number of prime divisors $\Gamma$ on $X$ with
\[
\adeg(\acherncl_1(\overline{H}_1) \cdots \acherncl_1(\overline{H}_1) \crest \Gamma) \leq h
\]
is finite

\item
For all $h'$, the number of prime divisors $\Gamma'$ on $X'$ with
\[
\adeg(\acherncl_1(\pi^*(\overline{H}_1)) \cdots \acherncl_1(\pi^*(\overline{H}_1)) \crest
\Gamma') 
\leq h'
\]
is finite.
\end{enumerate}
\end{Lemma}

\Proof
Let $X_0$ be the maximal Zariski open set of $X$ such that
$X_0$ is regular and $\pi$ is finite over $X_0$.
Then, $\codim(X \setminus X_0) \geq 2$.
We set $X'_0 = \pi^{-1}(X_0)$ and $\pi_0 = \rest{\pi}{X'_0}$.
Let $\Div(X)$ and $\Div(X')$ be the groups of Weil divisors on $X$ and
$X'$ respectively.
Then, a homomorphism $\pi^{\star} : \Div(X) \to \Div(X')$ is defined by
the compositions of homomorphisms:
\[
\Div(X) \to \Div(X_0) \overset{\pi_0^*}{\longrightarrow}
\Div(X'_0) \to \Div(X'),
\]
where $\Div(X) \to \Div(X_0)$ is the restriction map and
$\Div(X'_0) \to \Div(X')$ is defined by taking the Zariski closure of divisors.
Note that $\pi_*\pi^{\star}(D) = \deg(\pi)D$ for all $D \in \Div(X)$.

\medskip
First, we assume (1).
Note that the number of prime divisors in $X' \setminus X'_0$ is finite, so that
it is sufficient to show that
the number of prime divisors $\Gamma'$ on $X'$ with $\Gamma' \not\subseteq X' \setminus X'_0$
and
\[
\adeg(\acherncl_1(\pi^*(\overline{H}_1)) \cdots \acherncl_1(\pi^*(\overline{H}_1)) \crest
\Gamma') 
\leq h'
\]
is finite. By the projection formula,
\[
\adeg(\acherncl_1(\pi^*(\overline{H}_1)) \cdots \acherncl_1(\pi^*(\overline{H}_1)) \crest \Gamma')
= \adeg(\acherncl_1(\overline{H}_1) \cdots \acherncl_1(\overline{H}_1) \crest \pi_*(\Gamma')).
\]
Thus, by (1), the number of $(\pi_*(\Gamma'))_{\rm red}$ is finite.
On the other hand, the number of prime divisors
in $\pi^{-1}(\pi_*(\Gamma)_{\rm red})$ is finite. Hence we get (2).

\medskip
Next, we assume (2). Let $\Gamma$ be a prime divisor on $X$
with 
$\adeg(\acherncl_1(\overline{H}_1) \cdots \acherncl_1(\overline{H}_1) \crest \Gamma) \leq h$.
Then,
\[
\adeg(\acherncl_1(\pi^*(\overline{H}_1)) \cdots \acherncl_1(\pi^*(\overline{H}_1)) \crest
\pi^{\star}(\Gamma)) =
\deg(\pi) \adeg(\acherncl_1(\overline{H}_1) \cdots \acherncl_1(\overline{H}_1) \crest \Gamma)
\leq \deg(\pi) h.
\]
Thus, by (2), the number of $\pi^{\star}(\Gamma)$'s is finite.
Therefore, we get (1).
\QED

\bigskip
Let us go back to the proof of 
Proposition~\ref{prop:finite:divisor:fine:polarization}.
We use the notation in the above definition of strict finiteness.
By Lemma~\ref{lem:equiv:div:finite},
it is sufficient to show that the number of prime divisors $\Gamma'$ on $B'$ with
\[
\adeg(\acherncl_1(\pi^*(\overline{H}_1)) \cdots \acherncl_1(\pi^*(\overline{H}_d)) \crest \Gamma') \leq h
\]
is finite for all $h$.

There are $\QQ$-effective $C^{\infty}$-hermitian line bundles
$\overline{Q}_1, \ldots, \overline{Q}_d$ on $B'$ with
\[
\pi^*(\overline{H}_i) =  \mu^*(p_i^*(\overline{L}_i))
\otimes \overline{Q}_i
\]
for all $i$.
Note that
\begin{multline*}
\adeg(\acherncl_1(\pi^*(\overline{H}_1)) \cdots \acherncl_1(\pi^*(\overline{H}_d)) \crest \Gamma')
=
\adeg(\acherncl_1( \mu^*(p_1^*(\overline{L}_1))) \cdots \acherncl_1( \mu^*(p_d^*(\overline{L}_d))) \crest \Gamma') + \\
\sum_{i=1}^d
\adeg(\acherncl_1( \mu^*(p_1^*(\overline{L}_1))) \cdots
\acherncl_1( \mu^*(p_{i-1}^*(\overline{L}_{i-1}))) \cdot
\acherncl_1(\overline{Q}_i) \cdot
\acherncl_1(\pi^*(\overline{H}_{i+1})) \cdots
\acherncl_1(\pi^*(\overline{H}_d)) \crest \Gamma').
\end{multline*}
Moreover, since $\overline{Q}_i$ is $\QQ$-effective,
the number of prime divisors $\Gamma'$
with
\[
\adeg(\acherncl_1( \mu^*(p_1^*(\overline{L}_1))) \cdots
\acherncl_1( \mu^*(p_{i-1}^*(\overline{L}_{i-1}))) \cdot
\acherncl_1(\overline{Q}_i) \cdot
\acherncl_1(\pi^*(\overline{H}_{i+1})) \cdots
\acherncl_1(\pi^*(\overline{H}_d)) \crest \Gamma') < 0
\]
is finite for every $i$. Thus, we have
\[
\adeg(\acherncl_1(\pi^*(\overline{H}_1)) \cdots \acherncl_1(\pi^*(\overline{H}_d)) \crest \Gamma')
\geq
\adeg(\acherncl_1( \mu^*(p_1^*(\overline{L}_1))) \cdots \acherncl_1( \mu^*(p_d^*(\overline{L}_d))) \crest \Gamma')
\]
except finitely many $\Gamma'$.
On the other hand, by \cite[Proposition~4.1]{Mocycle},
the number of prime divisors $\Gamma''$ on $(\PP^1_{\ZZ})^d$ with
\[
\adeg(\acherncl_1(p_1^*(\overline{L}_1)) \cdots \acherncl_1(p_d^*(\overline{L}_d)) \crest \Gamma'')
\leq h
\]
is finite. Thus, we get our proposition.
\QED

\begin{Remark}
\label{rem:prop:finite:divisor:fine:polarization}
Let $X$ be a projective normal arithmetic variety of dimension $n$.
Let $\overline{H}_1, \ldots, \overline{H}_{n-2}$ be
nef $C^{\infty}$-hermitian line bundles on $X$ and $\overline{L}$
a $C^{\infty}$-hermitian line bundle on $X$.
If $\overline{L}$ is pseudo-effective, then we can expect
the number of prime divisors $\Gamma$ on $X$ with
\[
 \adeg(\acherncl_1(\overline{H}_1) \cdots \acherncl_1(\overline{H}_{n-2})
 \cdot \acherncl_1(\overline{L}) \crest \Gamma) < 0
\]
to be finite. If it is true, then 
Proposition~\ref{prop:finite:divisor:fine:polarization}
holds under the weaker assumption that
the polarization is fine.
\end{Remark}

\renewcommand{\theTheorem}{\arabic{section}.\arabic{Theorem}}
\renewcommand{\theClaim}{\arabic{section}.\arabic{Theorem}.\arabic{Claim}}
\renewcommand{\theequation}{\arabic{section}.\arabic{Theorem}.\arabic{Claim}}

\section{Height functions in terms of hermitian line bundles \\
with logarithmic singularities}

Let $K$ be a finitely generated field over $\QQ$ with
$d=\trdeg_{\QQ}(K)$.
Let $\overline{B} = (B; \overline{H}_1, \ldots, \overline{H}_d)$
be a fine polarization of $K$.
Let $X$ be a projective variety over $K$
and $L$ an ample line bundle on $X$.
Moreover, let $Y$  be a proper closed subset of $X$.
Let $(\mathcal{X}, \overline{\mathcal{L}})$ be a pair
of a projective arithmetic variety $\mathcal{X}$ and
a hermitian line bundle $\overline{\mathcal{L}}$ on
$\mathcal{X}$ with the following properties:
\begin{enumerate}
\renewcommand{\labelenumi}{(\arabic{enumi})}
\item
There is a morphism $f : \mathcal{X} \to B$ such that
the generic fiber of $f$ is $X$.

\item
$\mathcal{L}$ gives rise to $L$ on the generic fiber of $f$.

\item
$\mathcal{L}$ is ample with respect to the
morphism $f : \mathcal{X} \to B$.

\item
Let $\mathcal{Y}$ be a closed set of $\mathcal{X}$
such that $\mathcal{Y}$ gives rise to $Y$ on the generic fiber of
$\mathcal{X} \to B$.
Then the hermitian metric of $\overline{\mathcal{L}}$ has
logarithmic singularities along $\mathcal{Y}(\CC)$.
\end{enumerate}
For $x \in X(\overline{K}) \setminus Y(\overline{K})$, 
we denote by $\Delta_x$
the Zariski closure of the image of 
$\Spec(\overline{K}) \to X \to \mathcal{X}$.
The height of $x$ with respect to $\overline{\mathcal{L}}$
is defined by
\[
 h_{\overline{\mathcal{L}}}(x) =
\frac{\adeg(\acherncl_1(\rest{f^*(\overline{H}_1)}{\Delta_x})
\cdots \acherncl_1(\rest{f^*(\overline{H}_d)}{\Delta_x}) \cdot
\acherncl_1(\rest{\overline{\mathcal{L}}}{\Delta_x}))}{[K(x) : K]}.
\]
Note that since $\rest{\overline{\mathcal{L}}}{\Delta_x}$
has logarithmic singularities along $\mathcal{Y}(\CC) \cap
\Delta_x(\CC)$, the number
\[
\adeg(\acherncl_1(\rest{f^*(\overline{H}_1)}{\Delta_x})
\cdots \acherncl_1(\rest{f^*(\overline{H}_d)}{\Delta_x}) \cdot
\acherncl_1(\rest{\overline{\mathcal{L}}}{\Delta_x}))
\]
is well defined.
Then, we have the following proposition.

\begin{Proposition}
\label{prop:northcott:log:version}
\begin{enumerate}
\renewcommand{\labelenumi}{(\arabic{enumi})}
\item
Let us fix a positive integer $e$. Then, there is a constant $C$ such that
\[
 \# \{ x \in X(\overline{K}) \setminus Y(\overline{K}) \mid
 h_{\overline{\mathcal{L}}}(x) \leq h,\ [K(x) : K] \leq e \} \leq 
 C \cdot h^{d+1}
\]
for $h \gg 0$.

\item 
There is a constant $C'$ such that
$h_{\overline{L}}(x) \geq C'$ for all $x \in X(\overline{K}) \setminus
Y(\overline{K})$.
\end{enumerate}
\end{Proposition}

\Proof
We denote by $\Vert\cdot\Vert$ the hermitian metric of
$\overline{\mathcal{L}}$. 
Let $\overline{Q}$ be an ample $C^{\infty}$-hermitian line bundle on
$B$.
Then,
\[
 h_{\overline{\mathcal{L}} \otimes f^*(\overline{Q}^{\otimes n})}(x) =
 h_{\overline{\mathcal{L}}}(x) + 
 n \adeg(\acherncl_1(\overline{Q}) \cdot \acherncl_1(\overline{H}_1)
\cdots \acherncl_1(\overline{H}_d)).
\]
Thus, we may assume that $\mathcal{L}$ is ample on $\mathcal{X}$.
Moreover, replacing $\overline{\mathcal{L}}$ by 
$\overline{\mathcal{L}}^{\otimes n}$, we may assume that
$\mathcal{I}_{\mathcal{Y}} \otimes \mathcal{L}$ is
generated by global sections, where $\mathcal{I}_{\mathcal{Y}}$
is the defining ideal sheaf of $\mathcal{Y}$.
Let $s_1, \ldots, s_r$ be generators of
$H^0(\mathcal{X}, \mathcal{I}_{\mathcal{Y}} \otimes \mathcal{L})$.
We may view $s_1, \ldots, s_r$ as global sections
of $H^0(\mathcal{X}, \mathcal{L})$.
Then, 
$\mathcal{Y} = \{ x \in \mathcal{X} \mid s_1(x) = \cdots = s_r(x) = 0 \}$.
Here we choose a $C^{\infty}$-hermitian metric $\Vert\cdot\Vert_0$ of
$\mathcal{L}$ such that $\Vert s_i \Vert_0 < 1/e$ for all $i = 1, \ldots, r$.
We denote $(\mathcal{L}, \Vert\cdot\Vert_0)$ by $\overline{\mathcal{L}}^0$.

Here we claim
\[
 [K(x) : K]h_{\overline{\mathcal{L}}^0}(x) \geq
-\int_{\Delta_x(\CC)} \log \left(\max_i \{ \Vert s_i \Vert_0 \}\right)
c_1(f^*(\overline{H}_1))\wedge \cdots \wedge 
c_1(f^*(\overline{H}_d)).
\]
Indeed, we can find $s_j$ with $\rest{s_j}{\Delta_x} \not= 0$.
Thus,
\begin{multline*}
 [K(x) : K]h_{\overline{\mathcal{L}}^0}(x) =
\adeg( \acherncl_1(f^*(\overline{H}_1))
\cdots \acherncl_1(f^*(\overline{H}_d)) \crest
 \zeros(\rest{s_j}{\Delta_x})) \\
-\int_{\Delta_x(\CC)} \log \left(\Vert s_j \Vert_0 \right)
c_1(f^*(\overline{H}_1)) \wedge \cdots \wedge 
c_1(f^*(\overline{H}_d)).
\end{multline*}
Hence, we get our claim because
\[
 \adeg( \acherncl_1(f^*(\overline{H}_1))
\cdots \acherncl_1(f^*(\overline{H}_d)) \crest
 \zeros(\rest{s_j}{\Delta_x})) \geq 0
\quad\text{and}\quad
\Vert s_j \Vert_0 \leq \max_i \{ \Vert s_i \Vert_0 \}.
\]

Since $\Vert\cdot\Vert$ has logarithmic
singularities, if we set $g = \Vert\cdot\Vert/\Vert\cdot\Vert_0$,
then there is a positive constant $a, b$ such that
\[
 \vert \log(g) \vert \leq 
  a + b \log \left( - \log (\max_i \{
 \Vert s_i \Vert_0 \}) \right).
\]
Moreover,
\[
 \left\vert h_{\overline{\mathcal{L}}}(x) - h_{\overline{\mathcal{L}}^0}(x)
 \right\vert \leq \frac{1}{[K(x) : K]}
 \int_{\Delta_x(\CC)} \vert \log(g) \vert
c_1(f^*(\overline{H}_1)) \wedge \cdots \wedge 
c_1(f^*(\overline{H}_d)).
\]
Note that
\[
  \int_{\Delta_x(\CC)} 
c_1(f^*(\overline{H}_1)) \wedge \cdots \wedge 
c_1(f^*(\overline{H}_d)) = [K(x):K]\deg(\overline{B}),
\]
where ${\displaystyle \deg(\overline{B}) = \int_{B(\CC)}
c_1(\overline{H}_1) \wedge \cdots \wedge c_1(\overline{H}_d)}$
as in \S\S~\ref{subsec:polarization}.
Thus,
\[
 \frac{\left\vert h_{\overline{\mathcal{L}}}(x) -
 h_{\overline{\mathcal{L}}^0}(x)\right\vert}{\deg(\overline{B})}
 \leq a + 
b \int_{\Delta_x(\CC)}
\log \left( - \log (\max_i \{
 \Vert s_i \Vert_0 \}) \right) \frac{
c_1(f^*(\overline{H}_1)) \wedge \cdots \wedge 
c_1(f^*(\overline{H}_d))}{[K(x):K]\deg(\overline{B})}.
\]
On the other hand,
\begin{multline*}
  \int_{\Delta_x(\CC)}
\log \left( - \log (\max_i \{
 \Vert s_i \Vert_0 \}) \right) \frac{
c_1(f^*(\overline{H}_1)) \wedge \cdots \wedge 
c_1(f^*(\overline{H}_d))}{[K(x):K]\deg(\overline{B})} \\
\leq \log\left(
\int_{\Delta_x(\CC)}
- \log (\max_i \{
 \Vert s_i \Vert_0 \})\frac{
c_1(f^*(\overline{H}_1)) \wedge \cdots \wedge 
c_1(f^*(\overline{H}_d))}{[K(x):K]\deg(\overline{B})}
\right).
\end{multline*}
Hence, we obtain
\[
 \frac{\left\vert h_{\overline{\mathcal{L}}}(x) -
 h_{\overline{\mathcal{L}}^0}(x) \right\vert}{\deg(\overline{B})}
 \leq a + b \log\left(\frac{h_{\overline{\mathcal{L}}^0}(x)}{\deg(\overline{B})}\right).
\]
Note that there is a real number $t_0$ such that
$a + b \log(t) \leq t/2$ for all $t \geq t_0$.
Thus,
\[
 h_{\overline{\mathcal{L}}^0}(x) \leq \max \left\{ 
 \deg(\overline{B})t_0, 2 h_{\overline{\mathcal{L}}}(x)\right\}.
\]
Therefore, if $h \geq \deg(\overline{B})t_0/2$, then
$h_{\overline{\mathcal{L}}}(x) \leq h$ implies 
$h_{\overline{\mathcal{L}}^0}(x) \leq 2h$.
Hence, we get the first assertion by virtue of
\cite[Theorem~6.2.2]{Mocycle}.

Next let us see the second assertion. Since
\[
\Vert s_i \Vert = g \Vert s_i \Vert_0 \leq
\exp(a) \Vert s_i \Vert_0 \left(  
- \log (\max_j \{ \Vert s_j \Vert_0 \})\right)^b
\leq \exp(a) \Vert s_i \Vert_0 \left(  
- \log (\Vert s_i \Vert_0 )\right)^b
\]
and the function $t (-\log(t))^b$ is bounded above for $0 < t \leq 1$,
there is a constant $C$ such that $\Vert s_i \Vert \leq C$ for all $i$.
Thus, if we choose $s_i$ with $\rest{s_i}{\Delta_x} \not = 0$,
then
\begin{align*}
 [K(x) : K]h_{\overline{\mathcal{L}}}(x) & =
\adeg( \acherncl_1(f^*(\overline{H}_1))
\cdots \acherncl_1(f^*(\overline{H}_d)) \crest
 \zeros(\rest{s_i}{\Delta_x})) \\
& \qquad\qquad\qquad-\int_{\Delta_x(\CC)} \log \left(\Vert s_j \Vert \right)
c_1(f^*(\overline{H}_1)) \wedge \cdots \wedge 
c_1(f^*(\overline{H}_d)) \\
& \geq -\log(C)\int_{\Delta_x(\CC)}
c_1(f^*(\overline{H}_1)) \wedge \cdots \wedge 
c_1(f^*(\overline{H}_d))  \\
& = -\log(C) \deg(\overline{B})[K(x): K].
\end{align*}
Thus, we get (2).
\QED

\section{Faltings' modular height}
Let $K$ be a finitely generated field extension of $\QQ$ with 
$d = \trdeg_{\QQ}(K)$ and
$\overline{B} = (B; \overline{H}_1, \ldots, \overline{H}_d)$
a generically smooth polarization of $K$.
Let $A$ be a $g$-dimensional abelian variety over $K$.
Let $\hodge(A/K ; B)$ be the Hodge sheaf of $A$ with respect to $B$
(cf. \S\S~\ref{subsec:mod:sheaf:ab:var}).
Note that $\hodge(A/K;B)$ is invertible over $B_{\QQ}$
because $B_{\QQ}$ is smooth over $\QQ$.
Let $\Vert\cdot\Vert_{\rm Fal}$ be Faltings' metric
of $\hodge(A/K;B)$ over $B(\CC)$.
Here we set
\[
\overline{\hodge}^{\rm Fal}(A/K;B) = 
(\hodge(A/K;B), \Vert \cdot\Vert_{\rm Fal}),
\]
which is called {\em the metrized Hodge sheaf of $A$ with respect to $B$}.
By Lemma~\ref{lem:loc:L1:Faltings:metric}, 
the metric of $\overline{\hodge}^{\rm Fal}(A/K;B)$
is locally integrable. 
Thus, 
{\em the Faltings' modular height of $A$ with respect to the polarization 
$\overline{B}$} is defined by
\[
h_{\rm Fal}^{\overline{B}}(A) = \adeg\left(\acherncl_1(\overline{H}_1) \cdots \acherncl_1(\overline{H}_d)
\cdot \acherncl_1(\overline{\hodge}^{\rm Fal}(A/K;B)\right).
\]

\begin{Proposition}
\label{prop:comp:Faltings:hight:general}
Let $\pi : X' \to X$ be a generically finite morphism of
normal projective generically smooth
arithmetic varieties. Let $K$ and $K'$
be the function field of $X$ and $X'$ respectively.
Let $A$ be an abelian variety over $K$. Then,
there is an effective divisor $E$ on $X$ with the following properties:
\begin{enumerate}
\renewcommand{\labelenumi}{(\arabic{enumi})}
\item
$\pi_*\acherncl_1(\overline{\hodge}^{\rm Fal}(A\times_{K} \Spec(K')/K' ; X'))
=\deg(\pi)\acherncl_1(\overline{\hodge}^{\rm Fal}(A/K ; X)) + (E, 0)$.

\item
For a scheme $S$, we denote by $S^{(1)}$ the set of all codimension one points
of $S$.
Then,
\[
\{ x \in X^{(1)} \mid \text{$A$ has semi-abelian reduction at $x$} \}
\subseteq (X \setminus \Supp(E))^{(1)}.
\]
Moreover, if $A \times_K \Spec(K')$ has semi-abelian reduction
in codimension one, then
\[
\{ x \in X^{(1)} \mid \text{$A$ has semi-abelian reduction at $x$} \}
= (X \setminus \Supp(E))^{(1)}.
\]
\end{enumerate}
\end{Proposition}

\Proof
(1) Let $X_0$ be the maximal Zariski open set of $X$ such that
$X_0$ is regular and $\pi$ is finite over $X_0$.
Then, $\codim(X \setminus X_0) \geq 2$.
We set $X'_0 = \pi^{-1}(X_0)$ and $\pi_0 = \rest{\pi}{X'_0}$.
Let $\Div(X)$ and $\Div(X')$ be the groups of Weil divisors on $X$ and
$X'$ respectively.
Then, a homomorphism $\pi^{\star} : \Div(X) \to \Div(X')$ is defined by
the compositions of homomorphisms:
\[
\Div(X) \to \Div(X_0) \overset{\pi_0^*}{\longrightarrow}
\Div(X'_0) \to \Div(X'),
\]
where $\Div(X) \to \Div(X_0)$ is the restriction map and
$\Div(X'_0) \to \Div(X')$ is defined by taking the Zariski closure of divisors.
Note that $\pi_*\pi^{\star}(D) = \deg(\pi)D$ for all $D \in \Div(X)$.

Let $X_1$ (resp. $X'_1$) be a Zariski open sets of $X$ (resp. $X'$) such that
$\codim(X \setminus X_1) \geq 2$ (resp. $\codim(X' \setminus X'_1) \geq 2$) and
the N\'eron model $G$ (resp. $G'$) exists over $X_1$ (resp. $X'_1$).
Clearly we may assume that $X_1 \subseteq X_0$ and 
$\pi^{-1}(X_1) \subseteq X'_1$.
We set $X'_{2} = \pi^{-1}(X_1)$ and $G'_{2} = G' \times_{X'_1} X'_{2}$.
Since $G'_{2}$ is the N\'eron model of $A \times_K \Spec(K')$ over $X'_2$,
there is a homomorphism
$G \times_{X_1} X'_{2} \to G'_{2}$ over $X'_{2}$.
Thus, we get a homomorphism
\addtocounter{Claim}{1}
\begin{equation}
\label{eqn:prop:comp:Fal:hight:1}
\alpha : \pi^* \epsilon^* \left(\bigwedge^g \Omega_{G/X_1}\right) \to
{\epsilon'}^* \left(\bigwedge^g \Omega_{G'_{2}/X'_{2}}\right),
\end{equation}
where $\epsilon$ and $\epsilon'$ are the zero sections of $G$ and $G'$
respectively. 

Let $s$ be a non-zero rational section of $\hodge(A; X)$.
Then,
\[
\acherncl_1(\overline{\hodge}^{\rm Fal}(A/K ; X)) =
(\zeros(s), -\log \Vert s \Vert_{\rm Fal} ).
\]
Moreover, since $\pi^*(s)$ gives rise to a non-zero rational section of
$\hodge(A \times_{K} \Spec(K') ; X')$,
\[
\acherncl_1(\overline{\hodge}^{\rm Fal}(A\times_{K} \Spec(K')/K' ; X')) =
(\zeros(\pi^*(s)), -\pi^*(\log \Vert s \Vert_{\rm Fal}) ),
\]
where $\pi^*(\log \Vert s \Vert_{\rm Fal})$ is the pull-back of
$\log \Vert s \Vert_{\rm Fal}$ by $\pi$ as a function on a dense open set
of $X(\CC)$.
Let $\Gamma_1, \ldots, \Gamma_r$ be all prime divisors in
$X' \setminus X'_{2}$. Note that $\pi_*(\Gamma_i) = 0$ for
all $i$. Then, since 
\eqref{eqn:prop:comp:Fal:hight:1} is injective,
there is an effective divisor $E'$ and integers $a_1, \ldots, a_r$
such that
\[
\zeros(\pi^*(s)) = \pi^{\star}(\zeros(s)) + E' + \sum_{i=1}^r a_i \Gamma_i.
\]
Note that $E' = \sum_{x'} 
\length_{\OO_{X',x'}}(\Coker(\alpha)_{x'}) \overline{\{ x' \}}$,
where $x'$'s run over all codimension one points of $X'_2$.
Thus, since $\pi_*(\pi^{\star}(\zeros(s)), -\pi^*(\log \Vert s \Vert_{\rm Fal})) =
\deg(\pi) (\zeros(s), -\log \Vert s \Vert_{\rm Fal} )$, we have
\[
\pi_*\acherncl_1(\overline{\hodge}^{\rm Fal}(A\times_{K} \Spec(K')/K' ; X'))
=\deg(\pi)\acherncl_1(\overline{\hodge}^{\rm Fal}(A/K ; X)) + (\pi_*(E'), 0).
\]
Therefore, we get (1).

Next let us see (2).
We assume that $A$ has semi-abelian reduction at $x$.
Then, there is a open set $U$ such that
$x \in U$ and $\rest{G^{\rm o}}{U}$ is semi-abelian.
Thus, $\rest{G^{\rm o}}{U} \times_U \pi^{-1}(U)$ is semi-abelian.
Hence $\left(\rest{G'}{\pi^{-1}(U)}\right)^{\rm o}$ is isomorphic to
$\rest{G^{\rm o}}{U} \times_U \pi^{-1}(U)$. Thus
$x \not\in E_{\rm red}$.
Conversely, we assume that $A \times_{K} \Spec(K')$ has semi-abelian
reduction in codimension one and $x \not\in E_{\rm red}$.
Then, there is an open set $U \subset X_1$ such that $x \in U$ and
the homomorphism
\[
\alpha : \pi^* \epsilon^* \left(\bigwedge^g \Omega_{G/X_1}\right) \to
{\epsilon'}^* \left(\bigwedge^g \Omega_{G'_{2}/X'_{2}}\right)
\]
is an isomorphism over $\pi^{-1}(U)$, that is, so is
$\pi^* \epsilon^* \left(\Omega_{G/X_1}\right) \to
{\epsilon'}^* \left(\Omega_{G'_{2}/X'_{2}}\right)$
over $\pi^{-1}(U)$.
Thus, 
$G^{\rm o} \times_{X_1} X'_{2} \to (G'_{2})^{\rm o}$ is an isomorphism over 
$\pi^{-1}(U)$.
Therefore, $G^{\rm o}$ is semi-abelian over $U$.
\QED

\begin{Proposition}
\label{prop:isogeny:formula}
Let $\phi : A \to A'$ be an isogeny of abelian varieties
over $K$. Then
\begin{multline*}
\adeg\left(\acherncl_1(\overline{H}_1) \cdots \acherncl_1(\overline{H}_d)
\cdot \acherncl_1(\overline{\hodge}^{\rm Fal}(A'/K;B)\right)
- 
\adeg\left(\acherncl_1(\overline{H}_1) \cdots \acherncl_1(\overline{H}_d)
\cdot \acherncl_1(\overline{\hodge}^{\rm Fal}(A/K;B)\right) \\
= \frac{1}{2} \log(\deg(\phi)) \deg(\overline{B}) -
\adeg\left(\acherncl_1(\overline{H}_1) \cdots \acherncl_1(\overline{H}_d)
\crest D_{\phi} \right),
\end{multline*}
where $D_{\phi}$ is an effective divisor
given in \S\S~\rom{\ref{subsec:mod:sheaf:ab:var}} and
${\displaystyle \deg(\overline{B}) = 
\int_{B(\CC)} c_1(\overline{H}_1) \wedge \cdots
\wedge c_1(\overline{H}_d)}$ as in
\S\S~\rom{\ref{subsec:polarization}}.
\end{Proposition}

\Proof
This follows from the fact that 
$\overline{\hodge}^{\rm Fal}(A'/K;B) \otimes
(\OO_B(D_{\phi}), \deg(\phi)\vert\cdot\vert_{\rm can})$
is isometric to $\overline{\hodge}^{\rm Fal}(A/K;B)$.
\QED

\begin{Proposition}
\label{prop:height:dual}
If an abelian variety $A$ over $K$
has semi-abelian reduction in codimension one over $B$.
Then,
\[
\adeg\left(\acherncl_1(\overline{H}_1) \cdots \acherncl_1(\overline{H}_d)
\cdot \acherncl_1(\overline{\hodge}^{\rm Fal}(A/K;B)\right)
=
\adeg\left(\acherncl_1(\overline{H}_1) \cdots \acherncl_1(\overline{H}_d)
\cdot 
\acherncl_1(\overline{\hodge}^{\rm Fal}(A^{\vee}/K;B)\right),
\]
where $A^{\vee}$ is the dual abelian variety of $A$.
\end{Proposition}

\Proof
Let $\phi : A \to A^{\vee}$ be an isogeny over $K$ in terms
of ample line bundle on $A$.
Let $\phi^{\vee} : A \to A^{\vee}$ be the dual of $\phi$.
Then, by Proposition~\ref{prop:isogeny:formula},
\begin{multline*}
2 \left( \adeg\left(\acherncl_1(\overline{H}_1) \cdots
\acherncl_1(\overline{H}_d)
\cdot \acherncl_1(\overline{\hodge}^{\rm Fal}(A^{\vee}/K;B)\right)
- 
\adeg\left(\acherncl_1(\overline{H}_1) \cdots \acherncl_1(\overline{H}_d)
\cdot \acherncl_1(\overline{\hodge}^{\rm Fal}(A/K;B)\right) 
\right) \\
= \log(\deg(\phi))\deg(\overline{B}) -
\adeg\left(\acherncl_1(\overline{H}_1) \cdots \acherncl_1(\overline{H}_d)
\crest D_{\phi} + D_{\phi^{\vee}} \right).
\end{multline*}
On the other hand, by Lemma~\ref{lemma:mult:iso:isodual:deg},
$\mathcal{I}_{\phi} \cdot \mathcal{I}_{\phi^{\vee}} =
\deg(\phi) \OO_B$. Thus,
$(\OO_B(D_{\phi} + D_{\phi^{\vee}}), \vert\cdot\vert_{\rm can})$ is
isometric to $(\OO_B, \deg(\phi)^{-2}\vert\cdot\vert_{\rm can})$.
Therefore, we get our proposition.
\QED

\medskip
Let $K$ be a finitely generated field extension of $\QQ$ with 
$d = \trdeg_{\QQ}(K)$
and $\overline{B} = (B; \overline{H}_1, \ldots, \overline{H}_d)$
a polarization of $K$. 
Let $A$ be an abelian variety over a finite extension field $K'$
of $K$.
Let $m$ be a positive integer such that
$m$ has a decomposition $m=m_1 m_2$ with
$(m_1, m_2) = 1$ and $m_1, m_2 \geq 3$.
Let us consider a natural homomorphism
\[
 \rho(A,m) : \Gal(\overline{K}/K) \to \Aut(A[m](\overline{K}))
 \simeq \Aut((\ZZ/m\ZZ)^{2g}).
\]
Then, there is a Galois extension $K(A,m)$ of $K'$
with $\Ker \rho(A,m) = \Gal(\overline{K}/K(A,m))$.
Note that
\[
 \Gal(K(A,m)/K') = \Gal(\overline{K}/K)/\Ker \rho(A,m) \hookrightarrow
 \Aut((\ZZ/m\ZZ)^{2g}).
\]
Let $B''$ be a generically smooth, normal and projective arithmetic variety
with the following properties:
\begin{enumerate}
\renewcommand{\labelenumi}{(\roman{enumi})}
\item
The function field $K''$ of $B''$ is an extension of $K(A,m)$.

\item
The natural rational map $f : B'' \to B$ induced by $K \hookrightarrow K''$
is actually a morphism.
\end{enumerate}
Then, we have the following.

\begin{Proposition}
\label{prop:geom:fal:height}
\begin{enumerate}
\renewcommand{\labelenumi}{(\arabic{enumi})}
\item
The number
\[
\frac{1}{[K'': K]}
\adeg\left( \acherncl_1(\hodge(A \times_{K'} \Spec(K'')/K''; B'') \cdot
\acherncl_1(f^*(\overline{H}_1)) \cdots \acherncl_1(f^*(\overline{H}_1)) \right)
\]
does not depend on the choice of $m$ and $B''$, so that
we denote it by $h^{\overline{B}}_{\rm mod}(A)$.

\item
$h^{\overline{B}}_{\rm mod}(A) \leq h^{\overline{B}}_{\rm Fal}(A)$.
\end{enumerate}
\end{Proposition}

\Proof
These are consequences of 
Proposition~\ref{prop:semi:abelian:reduction},
Proposition~\ref{prop:comp:Faltings:hight:general}
and the projection formula.
\QED

\begin{Proposition}
\label{prop:height:prod}
Let $K$ be a finitely generated extension field of $\QQ$.
For abelian varieties $A$ and $A'$ over $K$,
$h^{\overline{B}}_{\rm Fal}(A \times_K A') = 
h^{\overline{B}}_{\rm Fal}(A) + h^{\overline{B}}_{\rm Fal}(A')$.
Moreover,
$h^{\overline{B}}_{\rm mod}(A \times_K A') = 
h^{\overline{B}}_{\rm mod}(A) + h^{\overline{B}}_{\rm mod}(A')$.
\end{Proposition}

\Proof
Let $\mathcal{A}$ and $\mathcal{A'}$ be the N\'eron models of $A$ and
$A'$ over $B_0$, where $B_0$ is a big open set of $B$.
Then, $\mathcal{A} \times_{B_0} \mathcal{A}'$ is the N\'eron model
of $A \times_K A'$ over $B_0$.
Thus,
\[
\acherncl_1(\overline{\hodge}^{\rm Fal}_{\mathcal{A} \times_{B_0}
\mathcal{A}'/B}) = 
\acherncl_1(\overline{\hodge}^{\rm Fal}_{\mathcal{A}/B_0}) +
\acherncl_1(\overline{\hodge}^{\rm Fal}_{\mathcal{A}'/B_0}).
\]
Hence, we get our lemma.
\QED

\section{Weak finiteness}
Let us fix positive integers $g$, $l$ and $m$ 
such that $m$ has a decomposition $m=m_1m_2$ with 
$(m_1, m_2) = 1$ and $m_1, m_2 \geq 3$.
Let $\abmod_{g,l,m,\QQ}$, $f : Y \to \abmod_{g,l,m}^*$,
$\overline{L}$, $n$ and $G \to Y$ be the same as 
in Proposition~\ref{thm:universal:extension:abel:semiable:with:level}.

Let $K$ be a finitely generated field extension of $\QQ$ with 
$d = \trdeg_{\QQ}(K)$ and let
$\overline{B} = (B; \overline{H}_1, \ldots, \overline{H}_d)$
be a generically smooth polarization of $K$.

Let $A$ be a $g$-dimensional and $l$-polarized abelian variety over a
finite extension $K'$ of $K$ with an $m$-level structure.
Let $x_A : \Spec(K') \to \abmod_{g,l,m}^*$ be the morphism
induced by $A$.
Moreover, let  $y_A : \Spec(K') \to \abmod_{g,l,m}^*  \times_{\ZZ} \Spec(K)$ 
be the morphism induced by $x_A$.
Let $\Delta_{A}$ be the closure of the image of $y_A$ in 
$\abmod^*_{g,l,m} \times _{\ZZ} B$.
Let $p : \abmod^*_{g,l,m} \times _{\ZZ} B \to \abmod^*_{g,l,m}$ and
$q : \abmod^*_{g,l,m} \times _{\ZZ} B \to B$ 
be the projections to the first factor and
the second factor respectively.
Here, we set
\[
h^{\overline{B}}_{\overline{L}}(A) =
\frac{1}{\deg(\Delta_A \to B)}
\adeg\left(
\acherncl_1(\rest{q^*(\overline{H}_1)}{\Delta_A}) \cdots 
\acherncl_1(\rest{q^*(\overline{H}_d)}{\Delta_A}) \cdot
\acherncl_1(\rest{p^*(\overline{L})}{\Delta_A})\right)
\]
which is nothing more than
the height of $y_A \in (\abmod_{g,l,m}^*  \times_{\ZZ} \Spec(K))(\overline{K})$
with respect to $\overline{L}$ and $\overline{B}$.
Then, we have the following proposition.

\begin{Proposition}
\label{prop:comp:fal:naive:height}
There is a constant $N(g,l,m)$ depending only on $g,l,m$
such that
\[
\vert h^{\overline{B}}_{\overline{L}}(A) - 
n h^{\overline{B}}_{\rm mod}(A) \vert
\leq \log(N(g,l,m)) \deg(\overline{B}).
\]
for every $g$-dimensional and $l$-polarized abelian variety $A$ over
$\overline{K}$ with an $m$-level structure, where
\[
 \deg(\overline{B}) = \int_{B(\CC)} c_1(\overline{H}_1) \wedge \cdots
 c_1(\overline{H}_d).
\]
\end{Proposition}

\Proof
Let $A$ be a $g$-dimensional and $l$-polarized abelian variety over
$\overline{K}$ with an $m$-level structure.
Let $K'$ be the minimal finite extension of $K$ such that
$A$, the polarization of $A$, the $m$-level structure
of $A$ are defined over $K'$.
Let $x_A : \Spec(K') \to \abmod_{g,l,m}^*$ be the morphism
induced by $A$. Moreover,
let $y_A : \Spec(K') \to \abmod_{g,l,m}^* \times_{\ZZ} B$
be the induced morphism by $x_A$.

Let $\Spec(K_1)$ be a closed  point of $Y \times_{\abmod_{g,l,m}^*} \Spec(K')$.
Then, we have the following commutative diagram:
\[
\xymatrix{
Y \ar[d]_{f} & & \Spec(K_1) \ar[ll] \ar[d] \\
\abmod^*_{g,l,m} & & \Spec(K') \ar[ll]_{x_A}
}
\]
Here, two $l$-polarized abelian varieties $A \times_{K'} \Spec(K_1)$ and 
$G \times_Y \Spec(K_1)$ with $m$-level structures 
gives rise to the same $K_1$-valued point of $\abmod^*_{g,l,m}$.
Thus, 
$A \times_{K'} \Spec(K_1)$ is isomorphic to $G \times_Y \Spec(K_1)$ over $K_1$
as $l$-polarized abelian varieties with $m$-level structures because
$m \geq 3$.
The above commutative diagram gives rise to the commutative diagram:
\[
\xymatrix{
Y \times_{\ZZ} B \ar[d] & & \Spec(K_1) \ar[ll] \ar[d] \\
\abmod_{g,l,m}^* \times_{\ZZ} B & & \Spec(K') \ar[ll]_{y_A}
}
\]
Let $B_1$ be a generic resolution of singularities of 
the normalization of $B$ in $K_1$. 
Note that a generic resolution of singularities
(a resolution of singularities over $\QQ$) exists by
Hironaka's theorem \cite{Hiro}.
Then, we have rational maps
$B_1 \dashrightarrow Y \times_{\ZZ} B$ and $B_1 \dashrightarrow \Delta_A$ such that
a composition
$B_1 \dashrightarrow \Delta_A \to \abmod_{g,l,m}^* \times_{\ZZ} B$ of rational maps
is equal to
$B_1 \dashrightarrow Y \times_{\ZZ} B \to \abmod_{g,m}^* \times_{\ZZ} B$.
Thus, there are a birational morphism $B_2 \to B_1$ of projective 
and generically smooth arithmetic varieties,
 a morphism $B_2 \to \Delta_A$ and a morphism $B_2 \to Y \times_{\ZZ} B$
with the following commutative diagram:
\[
\xymatrix{
B_1 \ar[d]_{\pi_1} & B_2 \ar[l] _{\gamma} \ar[rr]^{\beta} \ar[d]^{\alpha} & &
Y \times_{\ZZ} B \ar[d] _{f \times \operatorname{id}} \\
B & \Delta_A \ar[l] \ar[rr]^{\iota} &  & \abmod_{g,l,m}^* \times_{\ZZ} B
}
\]
Then,
\begin{align*}
h^{\overline{B}}_{\overline{L}}(A) & =
\frac{\adeg\left( \acherncl_1(\iota^*(p^*(\overline{L}))) \cdot
\acherncl_1(\iota^*(q^*(\overline{H}_1))) \cdots \acherncl_1(\iota^*(q^*(\overline{H}_1))) \right) }
{\deg(\Delta_A \to B)}
\\
& = \frac{\adeg\left( \acherncl_1(\alpha^*(\iota^*(p^*(\overline{L})))) \cdot
\acherncl_1(\alpha^*(\iota^*(q^*(\overline{H}_1)))) \cdots 
\acherncl_1(\alpha^*(\iota^*(q^*(\overline{H}_1)))) \right)}{\deg(B_2 \to B)} \\
& = \frac{\adeg\left( 
\acherncl_1(\beta^*((f \times \operatorname{id})^*(p^*(\overline{L}))))
\cdot
\acherncl_1(\gamma^*(\pi_1^*(\overline{H}_1))) \cdots 
\acherncl_1(\gamma^*(\pi_1^*(\overline{H}_1))) \right)}{\deg(B_2 \to B)}.
\end{align*}
On the other hand,
since $f^*(L) = \hodge_{G/Y}^{\otimes n}$ over $Y \times_{\ZZ} \Spec(\QQ)$,
there is an integer $N$ depending only on $g$, $l$ and $m$ such that
\[
N f^*(L) \subseteq \hodge_{G/Y}^{\otimes n} \subseteq (1/N) f^*(L) 
\]
on $Y$.
Thus,
\[
N \beta^*(f \times \operatorname{id})^*(L)
\subseteq (\hodge_{G \times_{\ZZ} B/Y
\times_{\ZZ} B})^{\otimes n}
\subseteq (1/N) \beta^*(f \times \operatorname{id})^*(L).
\]
Therefore,
\begin{multline*}
-\frac{\adeg(\acherncl_1(\gamma^*(\pi_1^*(\overline{H}_1))) \cdots 
\acherncl_1(\gamma^*(\pi_1^*(\overline{H}_1))) \crest (N))}{\deg(B_2 \to B)} +
h^{\overline{B}}_{\overline{L}}(A) \\
\leq \frac{n \adeg\left( \acherncl_1(\overline{\hodge}^{\rm Fal}_{G \times_Y B_2/B_2})) \cdot
\acherncl_1(\gamma^*({\pi_1}^*(\overline{H}_1))) \cdots 
\acherncl_1(\gamma^*({\pi_1}^*(\overline{H}_1))) \right)}{\deg(B_2 \to B)} \\
\leq \frac{\adeg(\acherncl_1(\gamma^*(\pi_1^*(\overline{H}_1))) \cdots 
\acherncl_1(\gamma^*(\pi_1^*(\overline{H}_1))) \crest (N))}{\deg(B_2 \to B)} +
h^{\overline{B}}_{\overline{L}}(A).
\end{multline*}
Note that
\[
\adeg(\acherncl_1(\gamma^*(\pi_1^*(\overline{H}_1))) \cdots 
\acherncl_1(\gamma^*(\pi_1^*(\overline{H}_1))) \crest (N)) =
\log(N)\deg(B_2 \to B) \deg(\overline{B}).
\]
By Proposition~\ref{prop:semi:abelian:reduction},
we can see that $A \times_{K'} \Spec(K_1)$ 
has semi-abelian reduction in codimension one
over $B_1$. On the other hand, by Proposition~\ref{prop:comp:Faltings:hight:general},
\[
\gamma_*(\acherncl_1(\overline{\hodge}^{\rm Fal}_{G \times_Y B_2/B_2})) =
\acherncl_1(\overline{\hodge}^{\rm Fal}(A \times_{K'} \Spec(K_1)/K_1; B_1)).
\]
Therefore, we get 
\[
\vert h^{\overline{B}}_{\overline{L}}(A) - 
n h^{\overline{B}}_{\rm mod}(A) \vert
\leq \log(N) \deg(\overline{B}).
\]
\QED

\begin{Corollary}
\label{cor:finite:K:value:point}
Let $K$ be a finitely generated field extension of $\QQ$ with $d = \trdeg_{\QQ}(K)$
and $\overline{B} = (B; \overline{H}_1, \ldots, \overline{H}_d)$ a 
generically smooth and fine polarization of $K$.
Let us fix a positive integer $l$.
Then, we have the following:
\begin{enumerate}
\renewcommand{\labelenumi}{(\arabic{enumi})}
\item
There is a constant $C$ such that
$C \leq h^{\overline{B}}_{\rm mod}(A)$
for any $l$-polarized abelian variety $A$ over $\overline{K}$.

\item
Let us fix a positive integer $e$. Then, there is a constant $C'$ such that
the number of the set 
\[
\left\{A \times_{K'} \Spec(\overline{K})  \left|
\begin{array}{l}
\text{$A$ is a $g$-dimensional and 
$l$-polarized abelian variety over} \\
\text{a finite extension $K'$ of $K$ with $[K' : K] \leq e$ and 
$h^{\overline{B}}_{\rm mod}(A) \leq h$.} 
\end{array}
\right\}\right/ \!\simeq_{\bar{K}}
\]
is less than or equal to $C' \cdot h^{d+1}$ for $h \gg 0$.
\end{enumerate}
\end{Corollary}

\Proof
Let us fix a positive number $m$ such that
$m$ has a decomposition $m=m_1m_2$ with $(m_1,m_2) = 1$
and $m_1,m_2 \geq 3$.
Then, any $l$-polarized abelian variety over $\overline{K}$ has
a $m$-level structure.
Thus, (1) is a consequence of Proposition~\ref{prop:northcott:log:version}
and Proposition~\ref{prop:comp:fal:naive:height}.

Let $A$ be an $l$-polarized abelian variety over a finite
extension $K'$ of $K$.
Let $K''$ be the minimal extension of $K'$ such that
$A[m](\overline{K}) \subseteq A(K'')$.
Then, $[K'':K'] \leq \#(\Aut(\ZZ/m\ZZ)^{2g})$.
Thus, by using Proposition~\ref{prop:northcott:log:version}
and Proposition~\ref{prop:comp:fal:naive:height},
we get (2). 
\QED

\section{Galois descent}
Let $A$ be a $g$-dimensional abelian variety over a field $k$.
Let $m$ be a positive integer prime to the characteristic of $k$.
Note that 
an $m$-level structure $\alpha$ of $A$ over a finite extension $k'$ of $k$
is an isomorphism 
$\alpha : (\ZZ/m\ZZ)^{2g} \to A[m](k')$.
If $k'$ is a finite Galois extension over $k$, then
we have a homomorphism
\[
 \epsilon(k'/k, A,\alpha) : \Gal(k'/k) \to \Aut((\ZZ/m\ZZ)^{2g})
\]
given by $\epsilon(k'/k, A,\alpha)(\sigma) = \alpha^{-1} \cdot \sigma_A \cdot
\alpha$, where
\[
 \begin{CD}
\sigma_A : A \times_{k} \Spec(k') 
@>{\operatorname{id}_A \times (\sigma^{-1})^a}>>
A \times_{k} \Spec(k')
 \end{CD}
\]
is the natural morphism arising from $\sigma$.
Note that $(\sigma\cdot\tau)_A = \sigma_A \cdot \tau_A$.

\begin{Lemma}
\label{lem:galois:descent}
Let $(A, \xi)$ and $(A',\xi')$ be polarized abelian varieties over a field $k$.
Let $m$ be a positive integer prime to the characteristic of $k$.
Let $\alpha$ and $\alpha'$ be $m$-level structures
of $A$ and $A'$ respectively over a finite Galois
extension $k'$ of $k$.
Let $\phi : (A,\xi) \times_k \Spec(k') \to (A',\xi') \times_k \Spec(k')$ be
an isomorphism as polarized abelian varieties over $k'$.
If $m \geq 3$, $\phi \cdot \alpha = \alpha'$ and
$\epsilon(k'/k, A, \alpha) = \epsilon(k'/k, A', \alpha')$,
then $\phi$ descents to an isomorphism
$(A,\xi) \to (A',\xi')$ over $k$.
\end{Lemma}

\Proof
For $\sigma \in \Gal(k'/k)$,
let us consider a morphism
\[
 \phi_{\sigma} = \sigma_{A'}^{-1} \cdot \phi \cdot \sigma_A
: A \times_k \Spec(k') \to A' \times_k \Spec(k').
\]
First of all, $\phi_{\sigma}$ is a morphism over $k'$.
We claim that $\phi_{\sigma} \cdot \alpha = \alpha'$.
Indeed, since $\alpha^{-1} \cdot \sigma_A \alpha = 
{\alpha'}^{-1} \cdot \sigma_{A'} \cdot \alpha'$, we have
\[
 \phi_{\sigma} \cdot \alpha = 
\sigma_{A'}^{-1} \cdot \phi \cdot \alpha \cdot \alpha^{-1} \cdot
\sigma_A \cdot \alpha 
= \sigma_{A'}^{-1} \cdot \alpha' \cdot {\alpha'}^{-1} \cdot
 \sigma_{A'} \cdot \alpha'
= \alpha'.
\]
Thus, $\phi_{\sigma}$ preserves the level structures of 
$A \times_k \Spec(k')$ and $A' \times_k \Spec(k')$.
Hence, since $m \geq 3$ and $\phi_{\sigma} \cdot \phi^{-1}$
preserve the polarization $\xi$ of $A$ over $k'$ (hence
$(\phi_{\sigma} \cdot \phi^{-1})^N = \operatorname{id}$ for $N \gg 1$), 
by virtue of Serre's theorem,
we have $\phi_{\sigma} = \phi$, that is,
\[
 \phi \cdot \sigma_A = \sigma_{A'} \cdot \phi
\]
for all $\sigma \in \Gal(k'/k)$.
Therefore, $\phi$ descents to an isomorphism
$(A,\xi) \to (A',\xi')$ over $k$.
\QED

\begin{Proposition}
\label{prop:finite:abel:var:alg:closure}
Let $B$ be an irreducible normal scheme such that $B$ is of finite type over $\ZZ$.
Let $K$ be the local ring at the generic point of $B$.
For a fixed $g$-dimensional polarized abelian variety $(C,\xi_C)$ 
over $\overline{K}$, we set
\[
\mathcal{S} = \left\{ (A,\xi) \left| 
\begin{array}{l} \text{$(A,\xi)$ is a polarized abelian variety
over $K$ with $(A,\xi) \times_{K} \Spec(\overline{K}) \simeq (C,\xi_C)$} \\
\text{and $A$ has semi-abelian reduction over $B$ in codimension one.}
\end{array} \right\}\right..
\]
Then, the number of isomorphism classes in $\mathcal{S}$ is finite.
\end{Proposition}

\Proof
For $(A,\xi)\in \mathcal{S}$,
let $B_A$ be a big open set of $B$ over which we have a semi-abelian
extension $\mathcal{X}_A \to B_A$ of $A$. 
Moreover, let $BR(A)$ be the set of
codimension one points $x$ of $B_A$ such that
the fiber of $\mathcal{X}_A$ over $x$ is not an abelian variety.

\begin{Claim}
For any $(A,\xi),(A',\xi') \in \mathcal{S}$,
$BR(A) = BR(A')$. 
\end{Claim}

Since $A \times_K \Spec(\overline{K}) \simeq 
A' \times_K \Spec(\overline{K}$),
there is a finite extension $K'$ of $K$ with
$A \times_K \Spec(K') \simeq A' \times_K \Spec(K')$.
Let $\pi : B' \to B$ be the normalization of $B$ in $K'$.
Then, $\mathcal{X}_A \times_{B_A} \pi^{-1}(B_A)$ is isomorphic
to $\mathcal{X}_{A'} \times_{B_{A'}} \pi^{-1}(B_{A'})$ over
$\pi^{-1}(B_A \cap B_{A'})$.
Thus, $\pi^{-1}(BR(A)) = \pi^{-1}(BR(A'))$.
Therefore, we get our claim.

\medskip
Let us fix a positive integer $m \geq 3$ and
$A_0 \in \mathcal{S}$. We set 
\[
 U = B \setminus \left( (B \times_{\ZZ} \Spec(\ZZ/m\ZZ) )\cup
       \Sing(B) \cup \bigcup_{x \in BR(A_0)} \overline{\{ x \}}\right).
\]
Then, $U$ is regular and of finite type over $\ZZ$.
The characteristic of the residue field of any point of $U$
is prime to $m$.
Moreover,
by the above claim, if we set $U_A = U \cap B_A$ for $A \in \mathcal{S}$,
then $\mathcal{X}_A$ is an abelian scheme over $U_A$ and
$\codim(U \setminus U_A) \geq 2$.

\begin{Claim}
There is a finite Galois extension $K'$ of $K$
such that for any $(A,\xi) \in \mathcal{S}$,
all $m$-torsion points of $A$ belong to $A(K')$. 
\end{Claim}

For $(A,\xi) \in \mathcal{S}$,
let $K_A$ be the finite extension of $K$ obtaining
by adding all $m$-torsion points of $A$ to $K$.
Let $V_A$ be the normalization of $U$ in $K_A$.
Then, it is well-known that
$V_A$ is \'etale over $U_A$.
Moreover, by virtue of the purity of branch loci
(cf. SGA~1, Expos\'e~X, Th\'er\`eme~3.1),
$V_A$ is \'etale over $U$.
Let $M$ be the union of finite extension $K'$ of $K$ such that
the normalization of $U$ in $K'$ is \'etale over $U$.
Then, it is easy to see that $M$ is a Galois extension of $K$.
Since $K_A \subseteq M$,
we have a continuous homomorphism
\[
 \rho_A : \Gal(M/K) \to \Aut(A[m](\overline{K})) \simeq
 \Aut((\ZZ/m\ZZ)^{2g})
\]
such that $\ker(\rho_A) = \Gal(M/K_A)$.
Since $\Gal(M/K) = \pi_1(U)$,
by \cite[Hermite-Minkowski theorem in Chapter~VI]{FalRat},
we have only finitely many continuous homomorphisms
\[
 \rho : \Gal(M/K) \to \Aut((\ZZ/m\ZZ)^{2g}).
\]
Thus, there are only finitely many Galois groups
$\{ \Gal(M/K_A) \}_{A \in \mathcal{S}}$.
Therefore, $\{ K_A \}_{A \in \mathcal{S}}$ is finite as a subfield
of $M$.
Thus, we get our claim.

\begin{Claim}
For any $(A,\xi), (A',\xi') \in \mathcal{S}$,
$(A,\xi) \times_K \Spec(K') \simeq (A',\xi') \times_K \Spec(K')$.
\end{Claim}

There is a finite Galois extension $K''$ of $K'$
such that an isomorphism
\[
\phi : (A,\xi) \times_K \Spec(K'') \to (A',\xi') \times_K \Spec(K'') 
\]
is given over $K''$.
Let $\alpha$ be an $m$-level structure of $A$ over $K''$ and
$\alpha' = \phi \cdot \alpha$.
Then, $\epsilon(K''/K', A \times_K \Spec(K'), \alpha) = 
\epsilon(K''/K', A' \times_K \Spec(K'), \alpha') = 1$ 
because all $m$-torsion points of $A$ and $A'$ are defined over $K'$.
Thus, $A \times_K \Spec(K'') \to A' \times_K \Spec(K'')$ descents to
an isomorphism
$(A,\xi) \times_K \Spec(K') \to (A',\xi') \times_K \Spec(K')$
by Lemma~\ref{lem:galois:descent}.

\medskip
Finally, let us see the number of isomorphism classes in $\mathcal{S}$ is finite.
Let us fix $(A_0,\xi_0) \in \mathcal{S}$ and
an $m$-level structure $\alpha_0$ of $A_0$ over $K'$.
Let $\phi_A : (A_0,\xi_0) \times_K \Spec(K') \to (A,\xi) \times_K \Spec(K')$ be
an isomorphism over $K'$. We set $\alpha_A = \phi_A \cdot \alpha_0$
and $\phi^{A}_{A'} = \phi_{A'} \cdot \phi_A^{-1} : A \times_K \Spec(K') \to
A' \times_K \Spec(K')$ for $(A,\xi), (A',\xi') \in \mathcal{S}$.
Then, $\alpha_{A'} = \phi^{A}_{A'} \cdot \alpha_A$.
Here let us consider a map
\[
 \gamma : \mathcal{S} \to \Hom(\Gal(K'/K), \Aut((\ZZ/m\ZZ)^{2g}))
\]
given by $\gamma(A) = \epsilon(K'/K, A, \alpha_A)$.
By Lemma~\ref{lem:galois:descent},
if $\gamma(A) = \gamma(A')$, then $(A,\xi) \simeq (A',\xi')$ over $K$.
Moreover, $\Hom(\Gal(K'/K), \Aut((\ZZ/m\ZZ)^{2g}))$ is a finite set.
Therefore, we get our proposition.
\QED

\section{Strong finiteness}
In this section, we give the proof of the main result of this note.

\begin{Theorem}
\label{thm:finite:prin:abel:bound:height}
Let $K$ be a finitely generated field over $\QQ$ with $d = \trdeg_{\QQ}(K)$.
Let $\overline{B} = (B; \overline{H}_1, \ldots, \overline{H}_d)$ be a 
generically smooth and strictly fine
polarization of $K$.
Then, for any numbers $c$,
the number of isomorphism classes of abelian varieties defined over $K$ 
with $h^{\overline{B}}_{\rm Fal}(A) \leq c$ is finite.
\end{Theorem}

\Proof
Let us consider the following two sets:
\begin{align*}
\mathcal{S}_0(c) & = \left\{ (A,\xi) \mid
\text{$(A,\xi)$ is a principally polarized abelian variety over $K$
with $h^{\overline{B}}_{\rm mod}(A) \leq 8c$}
\right\} \\
\mathcal{S}(c) & = \left\{ A \mid
\text{$A$ is an abelian variety over $K$
with $h^{\overline{B}}_{\rm Fal}(A) \leq c$}
\right\}
\end{align*}
Then, by Corollary~\ref{cor:finite:K:value:point},
$\left\{ (A,\xi) \times \Spec(\bar{K}) \mid 
\text{$(A,\xi) \in \mathcal{S}_0(c)$} \right\}/
\!\!\simeq_{\bar{K}}$
is finite.
By Zarhin's trick (cf. \cite[Expos\'e~VIII, Proposition~1]{Mordell}),
for an abelian variety $A$ over $K$,
$(A \times A^{\vee})^4$ is principally polarized. Moreover,
\[
h^{\overline{B}}_{\rm mod}((A \times A^{\vee})^4) =
8 h^{\overline{B}}_{\rm mod}(A).
\]
by Proposition~\ref{prop:height:dual} and 
Proposition~\ref{prop:height:prod}.
Thus, if $A \in \mathcal{S}(c)$, then
$(A \times A^{\vee})^4 \in \mathcal{S}_0(c)$.
Here, the number of isomorphism classes of
direct factors of $(A \times A^{\vee})^4 \times_K \Spec(\bar{K})$
is finite (cf. \cite[Expos\'e~VIII, Proposition~2]{Mordell}). Thus, 
$\{ A \times_K \Spec(\bar{K}) \mid A \in \mathcal{S}(c)\}/
\!\!\simeq_{{\bar{K}}}$ is
finite.
In particular, there is a constant $C$ such that
$C \leq h^{\overline{B}}_{\rm mod}(A)$ for all $A \in
\mathcal{S}(c)$.

Let $K_A$ be the minimal finite extension of $K$ such that
$A[12](\bar{K}) \subseteq A(K_A)$. Then,
$[K_A : K] \leq \#\Aut((\ZZ/12\ZZ)^{2g})$.
Let $B_A$ be a generic resolution of singularities of
the normalization of $B$ in $K_A$.
By Proposition~\ref{prop:semi:abelian:reduction}, $A \times_K \Spec(K_A)$ has
semi-abelian reduction in codimension one over $B_A$.
Thus, by Proposition~\ref{prop:comp:Faltings:hight:general},
there is an effective divisor $E_A$ on $B$ with
\[
h^{\overline{B}}_{\rm Fal}(A) - h^{\overline{B}}_{\rm mod}(A)
=\frac{\adeg(\acherncl_1(\overline{H}_1) \cdots \acherncl_1(\overline{H}_d) \crest E_A)}
{[K_A : K]}.
\]
Here $h^{\overline{B}}_{\rm mod}(A) \geq C$ for all $A \in
\mathcal{S}(c)$. Thus, we can find a constant $C'$ such that
\[
\adeg(\acherncl_1(\overline{H}_1) \cdots \acherncl_1(\overline{H}_d) \crest E_A)
\leq C'
\]
for all $A \in \mathcal{S}(c)$.
Therefore, 
by virtue of Proposition~\ref{prop:finite:divisor:fine:polarization},
there is a reduced effective divisor $D$ on $B$
such that, for all $A \in \mathcal{S}(c)$,
$A$ has semi-abelian reduction in codimension one over $B \setminus D$.
Hence, by Proposition~\ref{prop:finite:abel:var:alg:closure}, we have
our assertion.
\QED

\begin{Remark}
\label{rem:thm:finite:prin:abel:bound:height}
If the problem in Remark~\ref{rem:prop:finite:divisor:fine:polarization} 
is true,
then Theorem~\ref{thm:finite:prin:abel:bound:height} holds
even if the polarization $\overline{B}$ is fine.
\end{Remark}


\bigskip
\begin{thebibliography}{99}

\bibitem{Neron}
S. Bosch, W. L\"utkebohmert and M. Raynaud,
N\'eron models,
Springer, 1990.

\bibitem{FalRat}
G. Faltings , G. W\"ustholz et al.,
Rational points,
Aspect of Mathematics, Vol. E6, Vieweg, 1984.

\bibitem{FalChai}
G. Faltings and C. Chai,
Degeneration of abelian varieties,
Springer, 1990.

\bibitem{KMSemi}
S. Kawaguchi and A. Moriwaki,
Inequalities for semistable families for arithmetic varieties,
J. Math. Kyoto Univ. 36, 97--182 (2001).

\bibitem{Hiro}
H. Hironaka,
Resolution of singularities of an algebraic variety over a field of
characteristic zero, Ann. of Math. 79 (1964), I:109-203; II:205-326.

\bibitem{Moret}
L. Moret-Bailly,
Pinceaux de vari\'et\'es ab\'eliennes,
Ast\'erisque 129 (1985).

\bibitem{MoArht}
A. Moriwaki,
Arithmetic height functions over finitely generated fields,
Invent. math.,
140, 
(2000),
101--142.

\bibitem{Mocycle}
A. Moriwaki,
The number of algebraic cycles with bounded degree, preprint(math.AG/0209287).

\bibitem{Spgenus}
L. Szpiro, S\'eminaire sur les pinceaux de courbes de genre au moins deux,
ast\'erisque 86.

\bibitem{Mordell}
L. Szpiro, S\'eminaire sur les pinceaux arithm\'etiques: La conjecture de Mordell,
ast\'erisque 127.
\end{thebibliography}
\end{document}